\documentclass[final]{siamonline171218}
\pdfoutput=1
\usepackage{amsmath}
\usepackage{amssymb}
\usepackage{amsfonts}
\usepackage{graphicx}
\usepackage{color}

\newtheorem{exercise}{Exercise}
\newtheorem{exampl}{Example}
\newcommand{\comment}[1]{} 

\def\el{\vspace{11pt}\noindent} 

\def\bq{\begin{quotation}}
\def\eq{\end{quotation}}

\definecolor{darkgreen}{rgb}{0.2,0.5,0.2}

\def\d{\delta}

\def\VC{\Theta}

\def\s{\sigma}
\def\Si{\Sigma} 

\def\f{\varphi}

\def\F{\Phi}

\def\y{\psi}
\def\Y{\Psi}
\def\o{\omega}
\def\O{\Omega}

\newcommand{\V}[1]{ \mathbf{#1} }    

\newcommand{\Vo}{\V{0}}

\newcommand{\Vb}{\V{b}}
\newcommand{\Vc}{\V{c}}
\newcommand{\Vd}{\V{d}}

\newcommand{\Vm}{\V{m}}

\newcommand{\Vp}{\V{p}}

\newcommand{\Vr}{\V{r}}

\newcommand{\Vu}{\V{u}}
\newcommand{\Vv}{\V{v}}
\newcommand{\Vw}{\V{w}}
\newcommand{\Vx}{\V{x}}
\newcommand{\Vy}{\V{y}}
\newcommand{\Vz}{\V{z}}

\newcommand{\VGr}[1]{ \boldsymbol{#1} } %

\newcommand{\Vgd}{\VGr{\d}}

\newcommand{\Vgf}{\VGr{\f}}

\newcommand{\M}[1]{ \mathbf{#1} }  

\newcommand{\MA}{\M{A}}
\newcommand{\MB}{\M{B}}
\newcommand{\MC}{\M{C}}
\newcommand{\MD}{\M{D}}
\newcommand{\ME}{\M{E}}

\newcommand{\MI}{\M{I}}
\newcommand{\MJ}{\M{J}}
\newcommand{\MK}{\M{K}}
\newcommand{\ML}{\M{L}}

\newcommand{\MRr}{\M{R}}
\newcommand{\MS}{\M{S}}

\newcommand{\MU}{\M{U}}
\newcommand{\MV}{\M{V}}
\newcommand{\MW}{\M{W}}
\newcommand{\MX}{\M{X}}
\newcommand{\MY}{\M{Y}}
\newcommand{\MZ}{\M{Z}}

\newcommand{\MgVC}{\M{\VC}} 
\newcommand{\MgF}{\M{\F}}
\newcommand{\MgY}{\M{\Y}}
\newcommand{\MgO}{\M{\O}}

\newcommand{\MgSi}{\M{\Si}}

\newcommand{\cJ}{{\mathcal J}}
\newcommand{\cK}{{\mathcal K}}
\newcommand{\cL}{{\mathcal L}}

\newcommand{\Rn}[1]{\mathbb{R}^{#1}}
\newcommand{\Rmn}[2]{\mathbb{R}^{#1 \times #2}}

\def\Exp{\mathbb{E}}

\def\2nm#1{\|#1\|_2}

\def\Ra#1{\mathrm{Range}(#1)}

\def\Rk#1{\mathrm{Rank}(#1)}

\newcommand{\ars}[1]{\left[ \begin{array}{#1}}
\newcommand{\are}{\end{array} \right] }
\newcommand{\oars}[1]{\begin{array}{#1}}
\newcommand{\oare}{\end{array}}
\newcommand{\rars}[1]{\left( \begin{array}{#1}}
\newcommand{\rare}{\end{array} \right) }

\newcommand{\eqs}{\begin{eqnarray}}
\newcommand{\eqe}{\end{eqnarray}}
\newcommand{\eqsn}{\begin{eqnarray*}}
\newcommand{\eqen}{\end{eqnarray*}}

\def\defs{\begin{definition}}
\def\defe{\end{definition}}
\def\teos{\begin{theorem}}
\def\teoe{\end{theorem}}
\def\prfs{\begin{proof}}
\def\prfe{\end{proof}}
\def\exas{\begin{exampl}}
\def\exae{\end{exampl}}
\def\excs{\begin{exercise}}
\def\exce{\end{exercise}}
\def\cors{\begin{corollary}}
\def\core{\end{corollary}}
\def\lems{\begin{lemma}}
\def\leme{\end{lemma}}

\newcommand{\ens}{\begin{enumerate}}
\newcommand{\ene}{\end{enumerate}}

\newcommand{\its}{\begin{itemize}}
\newcommand{\ite}{\end{itemize}}

\newcommand{\des}{\begin{description}}
\newcommand{\dee}{\end{description}}

\def\wh{\widehat}

\def\wt{\widetilde}

\usepackage{algorithm,algpseudocode,float}
\newfloat{algorithm}{t}{lop}[section]

\usepackage{amssymb}
\usepackage[numbers]{natbib}
\usepackage{amsmath}
\usepackage{subeqnarray}
\usepackage{caption}
\usepackage[font={small,it}]{caption}
\usepackage{subcaption}
\usepackage{multirow,graphicx}
\usepackage[section]{placeins}

\newcommand{\mb}[1]{\mathbf{#1}}

\newcommand{\res}{\Vr} 
\newcommand{\Res}{\MRr} 

\newcommand{\randsrc}{\Vw}  
\newcommand{\Randsrc}{\MW}  
\newcommand{\Randdet}{\MV} 

\title{Randomized Approach to Nonlinear Inversion Combining 
Random and Optimized Simultaneous
Sources and Detectors\footnotemark[1]}

\author{Selin Aslan \footnotemark[2]\ ,
Eric de Sturler\footnotemark[2]\ ,
and Misha E. Kilmer\footnotemark[3]}

\begin{document}
\maketitle
\renewcommand{\thefootnote}{\fnsymbol{footnote}}
\footnotetext[1]{This material is based upon work supported by the National Science
Foundation under Grants No. {NSF DMS} 1217156 and 1217161,
and {NIH} R01-CA154774.}
\footnotetext[2]{Department of Mathematics, Virginia Tech, Blacksburg, VA 24061.}
\footnotetext[3]{Department of Mathematics, Tufts University, Medford, MA 02115.}
\renewcommand{\thefootnote}{\arabic{footnote}}

\begin{abstract}

In partial differential equations-based (PDE-based) inverse problems with
many measurements, many large-scale discretized PDEs must be solved for
each evaluation of the misfit or objective function.
In the nonlinear case, evaluating the
Jacobian requires solving an additional set of systems.
This leads to a tremendous computational
cost, and this is by far the dominant cost for these problems.
Several authors have proposed randomization and stochastic
programming techniques to drastically reduce the number of
system solves by estimating the objective function using only
a few appropriately chosen random linear combinations of the sources.
While some have reported good solution quality at a greatly reduced
cost, for our problem of interest, diffuse
optical tomography, the approach often
does not lead to sufficiently accurate solutions.

We propose two improvements.
First, to efficiently exploit Newton-type methods, we modify the 
stochastic estimates to include random linear combinations of detectors,
drastically reducing the number of adjoint solves.
Second, after solving to a modest tolerance,
we compute a few simultaneous sources and detectors 
that maximize the Frobenius norm of the sampled Jacobian
to improve the rate of convergence and obtain more
accurate solutions. We complement these
optimized simultaneous sources and detectors by
random simultaneous sources and detectors
constrained to a complementary subspace.
Our approach leads to solutions of the same quality as obtained using all
sources and detectors but at a greatly reduced computational cost,
as the number of large-scale linear systems to be
solved is significantly reduced.

\end{abstract}

\begin{keywords}
DOT, PaLS, stochastic programming, randomization, inverse problems, optimization
\end{keywords}

\begin{AMS}
65F22, 65N21, 65N22, 65M32, 62L20, 90C15
\end{AMS}

\section{Introduction} The solution of nonlinear inverse problems requires solving many
large-scale discretized PDEs in the evaluation of the forward problem. In
parameterized inverse problems, we can compute the response of the system
for a particular input by numerically solving the PDE. The forward model used in this
paper, see Section~\ref{Sec:SimultOpt_backg}, is already regularized using the parametric level set (PaLS) approach
\citep{Aghasi_etal11}, and we focus on
efficiently solving the nonlinear least squares problem
\eqs \label{eq:NonlinLS}
  \min\limits_{\Vp} f(\Vp) := \min\limits_{\Vp} \frac{1}{2}\|\mathbb{M}(\Vp) - \V{d} \|_{2}^2,
\eqe
where $\mathbb{M}(\Vp)$ is the vector of computed measurements given
by the
forward model for the parameter vector $\Vp$, and $\Vd$ is the vector of measured
data at the detectors.

Each evaluation of $f(\Vp)$ requires the solution of the PDE for
all inputs and each frequency. Moreover, to efficiently compute derivative information using the co-state approach \cite{Vogel2002}, we also need to solve linear systems with the
adjoint for each detector and each frequency. This leads to an enormous
computational bottleneck, as rapid advances in technology allow for large
numbers of sources and detectors. Multiply this by the number of frequencies,
and the number of linear systems to solve in the solution of (\ref{eq:NonlinLS}) is
very large indeed. For the main application discussed in this paper, diffuse optical tomography (DOT), the number of sources and the number of detectors may each be a thousand or more; the number of frequencies used is typically modest (less than ten) \cite{StuGuKilChatBeatOCon}.

To solve the minimization problem
(\ref{eq:NonlinLS}), we use the Trust region algorithm with Gauss-Newton REGularized model
Solution (TREGS) \citep{StuKil11c} that has proven very effective for parameterized problems of
the type we consider in this paper. In \cite{StuGuKilChatBeatOCon}, we use reduced order models
to approximate both the function evaluation as well as its derivatives
to compute regularized Gauss-Newton steps in TREGS. Here,
we explore an alternative approach, following the work
by Haber, Chung, and Herrmann \cite{HabeChunHerr2012}.
The main idea in their paper was to drastically reduce the number of
systems to be solved by exploiting randomization \cite{HabeChunHerr2012},
posing the problem as a stochastic optimization problem \cite{ShapDentRusz2009}. In their
approach, the misfit or objective function is estimated using only
a few
random linear combinations of the sources,
referred to as \emph{random simultaneous sources}, that are kept fixed over many optimization steps. In \cite{ShapDentRusz2009}, this approach is referred to as the Sample Average Approximation (SAA) method.

The use of random simultaneous sources has been well-studied in several papers; see \cite{Beasley, NeelamaniKrohnKrebRombDeffAnd, KrebsAndersHinkNeelLeeBaumLac, FarbodDoelAscher2, LeeuwenArawkHerr} and the references therein. While replacing the original objective function by the
stochastic optimization problem
seems to work well for direct current resistivity and seismic tomography \cite{HabeChunHerr2012},
we find that the approach does not lead to accurate recovery of the parameters
for the DOT problem. Therefore, we propose {\it two innovations} to the
use of random simultaneous sources.

First, we extend the idea of random simultaneous sources to
the randomized treatment of the detector solves for
efficiently computing the Jacobian in Newton-type methods. Second, we propose
to combine random simultaneous sources and detectors with optimized simultaneous
sources and detectors to best capture the sensitivity and hence obtain
more accurate estimates of the dominant singular components of the
Jacobian and the corresponding components of the gradient.

The first innovation drastically reduces the cost of
Newton-type methods.
In particular, we derive a stochastic optimization problem,
analogous to randomized simultaneous sources, that allows us
to reduce the number of adjoint solves for the detectors.

The second innovation avoids stagnation in the residual norm
decrease of the stochastic optimization approach due to poor
or less effective estimates of derivative information.
This is typically more important closer to the solution than
early in the optimization, and several authors have addressed this
problem by dynamically varying the sample size in the stochastic algorithm
or increasing it slowly; see, for example, \cite{ByrdChinNocedalWu,FarbodDoelAscher2}.
We propose an alternative method to improve the estimates of derivative information
that keeps the sample size fixed (and small) for efficiency. Comparing the
two approaches in detail is future work.
For the DOT problem, using random simultaneous sources and detectors
does provide moderately accurate parameter solution estimates
at a drastically reduced number of linear system solves.
Thus, in our new approach, we first solve with a fixed set of random simultaneous sources and detectors to an intermediate tolerance. After reaching this intermediate tolerance, we compute a
few simultaneous sources and detectors that maximize the Frobenius norm of the
sampled Jacobian (see Section \ref{Sec:SAA_Opt});
we refer to these as \emph{optimized simultaneous sources and detectors}.
We complement these optimized sources and detectors by random simultaneous sources and
detectors constrained to a complementary subspace (see Section \ref{Sec:SAA_Opt}).
After this update, the optimization converges rapidly to a solution of the same quality as
obtained using all sources and detectors.
Our use of optimized simultaneous sources and detectors
is based on two motivations. First, the regularized model problem solves in
TREGS \cite{StuKil11c} focus on the directions corresponding to the large singular values
of the Jacobian. Second, the directions corresponding to the large singular
values are best informed by the
data. More details follow at the end of Section \ref{Sec:SimultOpt_backg}.

This paper is organized as follows. In Section \ref{Sec:SimultOpt_backg}, we briefly review DOT, PaLS, and TREGS. In Section \ref{Sec:SAA_Opt}, we introduce an alternative stochastic problem that includes random
simultaneous detectors, to reduce the number of adjoint solves. In Section~\ref{sec:OptSrcDet_new},
we introduce optimized simultaneous sources and detectors combined with random simultaneous
sources and detectors constrained to a complementary subspace. We also give an outline of our implementation strategies. In Section \ref{Sec:SimultOpt_NumExp}, we demonstrate the effectiveness of combining random and optimized simultaneous sources and detectors using a 2D and a 3D experiment. Finally, we draw some conclusions and discuss future work in Section \ref{Sec:SimultOpt_con}.

\section{Background on DOT, PaLS, and TREGS}\label{Sec:SimultOpt_backg} We assume that the region to be imaged is a rectangular prism with
sources and detectors on the top and or the bottom.
%
We consider the diffusion model for the photon flux $\eta(\mb{x})$ obtained by an input source $g(\mb{x})$ as in \cite{arridge}. Let the diffusion (or the scattering) and the absorption coefficients be given by $D(\mb{x})$ and $\mu(\mb{x})$, respectively. Then, the mathematical model of the problem in the frequency domain is given by
\begin{align} \label{eq: PDE}
-\nabla \cdot (D(\Vx) & \nabla\eta(\Vx))+ \mu(\Vx)\eta(\Vx)+ \frac{\imath\omega}{\nu}\eta(\Vx)= g(\mb{x}), &\\  \nonumber
&\text{for} \quad \Vx=(x_1,x_2,x_3)^T\text{and} -a<x_1<a,\:\:\:-b<x_2<b, \:\:\:0<x_3<c, &\\  \nonumber
\eta(\Vx)=0\:\:  & \text{if}\:\: 0\leq x_3 \leq c\:\:\: \text{and} \:\:\: \text{either}\:\:\: x_1=\pm a, \:\:\:\text{or}\:\:\:\: x_2=\pm b, &\\  \nonumber
0.25\eta(\Vx)+&\frac{D(\Vx)}{2}\frac{\partial \eta(\Vx)}{\partial\xi}= 0  \:\:  \text{for}\:\:\:x_3=0,\:\:or\:\: x_3=c,
\end{align}
where $\xi$ is the outward unit normal,  $\omega$ is the frequency modulation of light,
and $\nu$ is the speed of light in the medium.

Assuming that the diffusion coefficient is known (a common assumption for breast imaging), we use measurements and the forward model to recover the absorption coefficient of the medium, which can be used to distinguish healthy tissue from tumors \cite{boas2001imaging}. Typical inversion methods would optimize for the desired physical quantity over a collection of grid points/voxels resulting in a parameter vector with at least $O(10^6)$ unknowns. Instead, we assume that the absorption field, $\mu(\mb{x})$, is expressible as $\mu(\mb{x};\Vp)$
with a modest number of (unknown) parameters, $\Vp=[p_1, p_2, \ldots, p_{n_p}]^T$, where $n_p$ is the number of parameters. We use the PaLS approach \citep{Aghasi_etal11, StuGuKilChatBeatOCon} and parameterize the absorption $\mu(\mb{x};\Vp)$ as follows.

Let $\varphi:\mathbb{R}^+\rightarrow \mathbb{R}$ be a smooth, compactly supported radial basis function (CSRBF)\footnote{The CSRBF used here, $\varphi(r)$, is the Wendland function
$\y_{2,1}(r) =(1-r)^4(4r+1)$ \cite[Table 1]{Aghasi_etal11}.}, $\gamma$ be a positive, small, real number, and $\|\mb{x}\|^\dagger:=\sqrt{\|\mb{x}\|_2^2+\gamma^2}$ denote the (regularized) Euclidean norm. Then the PaLS function $\phi$ with a vector of unknown parameters $\Vp$ consisting of expansion coefficients $\alpha_j$, dilation coefficients $\beta_j$, and center locations $\pmb{\chi}_j$ is defined as

\begin{eqnarray}
\phi(\mb{x},\mb{p}):=\sum\limits_{j=1}^{m_0}\alpha_j\varphi(\|\beta_j(\mb{x}-\pmb{\chi}_j)\|^\dagger).
\end{eqnarray}
The PaLS approach uses an approximate Heaviside function
$H_\epsilon(r)$, where r is a scalar, to create a differentiable, but sharp transition from anomaly to background. The absorption $\mu(\mb{x},\mb{p})$ takes the value $\mu_{in}(\mb{x})$ if $\mb{x}$ is inside the region and  $\mu_{out}(\mb{x})$ if  $\mb{x}$ is outside the region,

\begin{equation}\label{absorption const}
\mu(\mb{x},\mb{p})=\mu_{in}(\mb{x})H_\epsilon(\phi(\mb{x},\mb{p})-c)+\mu_{out}(\mb{x})(1-H_\epsilon(\phi(\mb{x},\mb{p})-c)),
\end{equation}
where $c \in \mathbb{R}$ is a chosen cut-off parameter for the level set.

\begin{figure}
    \centering
    \begin{subfigure}[b]{0.32\textwidth}
        \includegraphics[width=0.95\textwidth]{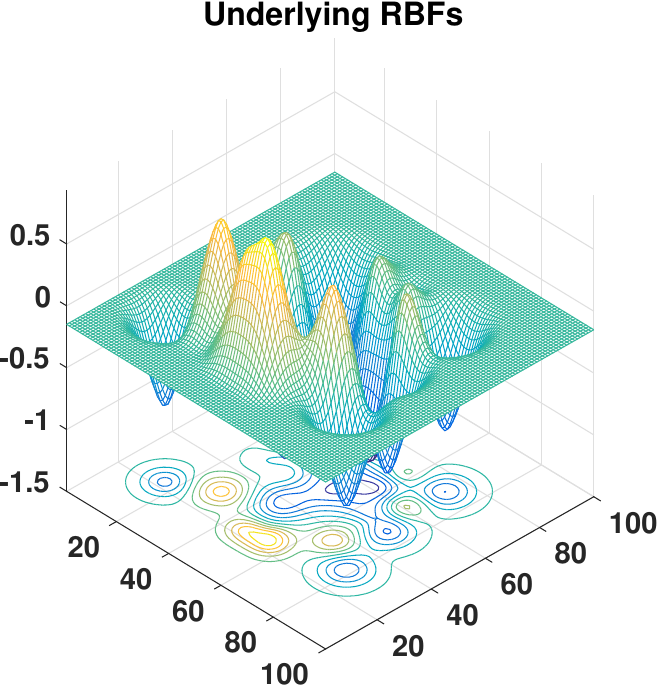}
        \subcaption{}
    \end{subfigure}
    \hspace{1in}
    \begin{subfigure}[b]{0.32\textwidth}
        \includegraphics[width=0.95\textwidth]{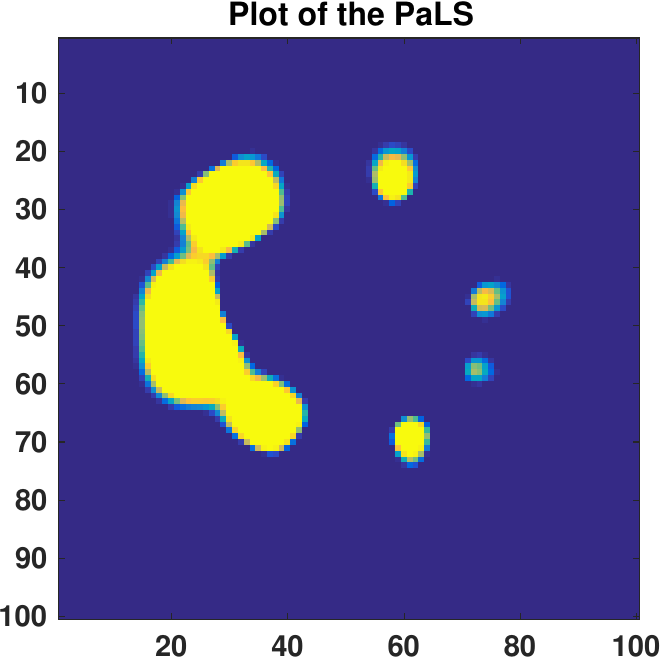}
        \subcaption{}
    \end{subfigure}
  \caption{(a) Surface and contour plot of a test anomaly on
  $100 \times 100$ mesh with 25 basis functions where the cut off is at $c=0.15$.
  (b) The PaLS function of the test anomaly on the left. If
  $\phi(\mb{x},\mb{p})\geq 0.15$, then $\mb{x}$ is inside the anomaly
  (light) and if $\phi(\mb{x},\mb{p})<0.15$, then $\mb{x}$ is
  outside the anomaly (dark). }
  \label{fig:Pals}
\end{figure}

Figure \ref{fig:Pals} illustrates how PaLS represents the absorption field. Using PaLS, edges and complex boundaries can be captured with relatively few basis functions. Moreover, the PaLS representation with a modest
number of basis functions regularizes the problem, hence no further regularization is needed.
Since there is no point to reduce the misfit below the
(known or estimated) norm of the
noise in the data, we stop the optimization when the objective function
reaches this noise level. This is called the discrepancy principle \cite{Vogel2002}.
For further discussion of the PaLS parameters for DOT, we refer the reader to \citep{Aghasi_etal11, StuGuKilChatBeatOCon}.

Let $n_d$, $n_s$, and $n_\omega$ denote the number of detectors, sources, and frequencies, respectively. The discretization of (\ref{eq: PDE}) leads to computed measurements, $\Vm_i(\omega_j,\Vp) \in \mathbb{C}^{n_{d}}$, for each source term, $\Vb_i$,
\begin{equation}\label{Eq:compmeas}
\Vm_i(\omega_j,\Vp)  =  \MC^T\left( \frac{\imath\omega_j}{\nu} \ME +\MA(\Vp)\right)^{-1} \Vb_i,
\end{equation}
where the rows of $\MC^T$ correspond to the detectors\footnotemark \footnotetext{In practice, we also split $\Vm_i$ in its real and imaginary parts.}. $\MA(\Vp)$ derives from a finite difference discretization of the diffusion and absorption terms in (\ref{eq: PDE}), and $\ME$ derives from the frequency term in (\ref{eq: PDE}).
$\ME$ is almost the identity except that it has zero rows for points on the boundary, $x_3=0,\: x_3=c,$ in (\ref{eq: PDE});
so, $\ME$ is singular.

For simplicity, we consider the nonlinear
residual for a single frequency, $\omega_j=0$.
In vector form, the residual is defined as follows
\begin{equation}\label{NonResidual}
 \Vr(\Vp)=  \begin{bmatrix}
          \Vr_1(\Vp) \\
           \vdots \\
           \Vr_{n_{s}}(\Vp) \\
             \end{bmatrix}= \begin{bmatrix}
 \Vm_1(\Vp)- \Vd_1\\
          \vdots \\
 \Vm_{n_s}(\Vp)- \Vd_{n_{s}}
  \end{bmatrix}=
 \begin{bmatrix}
 \MC^T \MA^{-1}(\Vp)\Vb_1 - \Vd_1\\
          \vdots \\
          \MC^T \MA^{-1}(\Vp)\Vb_{n_{s}} - \Vd_{n_{s}}
  \end{bmatrix},
\end{equation}
where $\Vr_i \in \mathbb{R}^{n_{d}}$, $\Vd_i$ is the data vector
with the measurements from the detectors corresponding to
source $\Vb_i$,
and the nonlinear least squares problem (\ref{eq:NonlinLS}) becomes
\begin{equation}\label{NonLeastSq1}
  \min_{\Vp} \frac{1}{2} \| \Vr(\Vp) \|_2^2.
\end{equation}
Let  $\MJ$ be the Jacobian of $\Vr(\Vp)$,
\begin{equation}\label{eq:Jacobian}
\MJ=\frac{\partial \Vr(\Vp)}{\partial \Vp} =
    \ars{ccc}
     \displaystyle \frac{\partial \Vr(\Vp)}{\partial \Vp_1} & \ldots &  \displaystyle \frac{\partial \Vr(\Vp)}{\partial \Vp_{n_{p}}}
    \are,
\end{equation}
where the components of $\MJ$ are given by the small vectors
\begin{equation}\label{eq:Jacobiancomp}
\MJ_{jk}(\Vp)=\frac{\partial}{\partial \Vp_k} (\MC^T \MA^{-1}(\Vp)\Vb_j)=-\MC^T \MA^{-1}(\Vp)\frac{\partial \MA(\Vp)}{\partial \Vp_k}\MA^{-1}(\Vp)\Vb_j \in \mathbb{R}^{n_d} .
\end{equation}
Evaluating the objective function at $\Vp$ requires solving $n_s\cdot n_\omega$ large linear systems. Once $\Vr(\Vp)$ and $\MA^{-1}(\Vp)\Vb_j$ are available, evaluating $\MJ$ using the co-state approach~\cite{Vogel2002} requires solving an additional $n_{d} \cdot n_{\omega}$  adjoint 
systems for the detectors. As a result, standard optimization approaches require $O(10^3-10^4)$ large linear system solves at each optimization step. The size of a realistic linear system is at least $O(10^6)$. This leads to an enormous computational bottleneck, and new computational techniques are needed.

We use TREGS \cite{StuKil11c} to solve the nonlinear least
squares problem (\ref{NonLeastSq1}). The TREGS algorithm combines a
trust region method with a regularized minimization of the Gauss-Newton (GN)
model \cite{dennis1996numerical}.  The local (GN) model at the
current parameter vector, $\Vp_c$, is given by
\begin{equation}\label{GN}
f(\Vp_c+\pmb{\d}) \approx m_{GN}(\Vp_c+\pmb{\d})= \frac{1}{2}
  {\Vr_c}^T\Vr_c+\Vr_c^T\MJ_c\pmb{\d}+\frac{1}{2} \pmb{\d}^T{\MJ_c}^T\MJ_c\pmb{\d},
\end{equation}
and its minimization is equivalent to the least squares
problem
\begin{equation}\label{Eq:NLLS}
  \min\limits_{\pmb{\d}} \|\MJ_c\pmb{\d}+\Vr(\Vp_c)\|_{2}^2.
\end{equation}
The TREGS algorithm favors updates corresponding to
(1) the large singular values and (2) the left singular vectors with
large components in $\Vr$ as determined by a generalized cross
validation-like (GCV) criterion. Since the Jacobian tends to be
ill-conditioned, the emphasis on large singular values
leads to relatively small steps that provide relatively
large reductions in the GN model (\ref{GN}).
We refer the reader to  \cite{StuKil11c} for more details of TREGS. 
\section{A Randomized Approach}\label{Sec:SAA_Opt}
We recast the nonlinear least squares problem as a stochastic optimization problem using randomization to drastically reduce the number of
large linear system solves in (\ref{NonResidual}) and (\ref{eq:Jacobiancomp}). The columns of $\MB=[\Vb_1,\cdots,\Vb_{n_{s}}]$ are source terms, and we refer to any linear combination of these sources as a simultaneous source.
Simultaneous random sources, $\MB \randsrc$, with $\randsrc \in \mathbb{R}^{n_s}$ a random vector, have been used in several areas \cite{Beasley, MortonOber, NeelamaniKrohnKrebRombDeffAnd,HabeChunHerr2012}.
In this section, we introduce the concept of \emph{optimized simultaneous sources and
detectors} to improve the rate of convergence of the optimization and the quality of the inverse solution.

\subsection{A Stochastic Optimization Approach}
To recast (\ref{NonResidual})--(\ref{NonLeastSq1}) as a stochastic
optimization problem, we first write the residual in matrix form. For a single frequency, we get
\begin{eqnarray}\label{MartixR}
\MRr(\Vp)=[ \Vr_1(\Vp)  \:\: \: \Vr_2(\Vp)\:\: \cdots \:\:\: \Vr_{n_{s}}(\Vp) ]= \MC^T \MA^{-1}(\Vp)\MB - \MD,
\end{eqnarray}
where
the vectors $\Vr_i \in \mathbb{R}^{n_{d}}$ are defined in (\ref{NonResidual}),
and consequently $\Vr(\Vp)=\text{vec}(\MRr(\Vp))$.\footnotemark
\footnotetext{For multiple frequencies, we need to compute the residual for each frequency, $[\MRr(\omega_1, \Vp) \  \MRr(\omega_2, \Vp)\:\: \cdots]=[\MC^T \MA^{-1}(\omega_1, \Vp) \: \MB - \MD_1 \ \ \ \  \MC^T \MA^{-1}(\omega_2, \Vp) \: \MB - \MD_2 \:\:\:\: \cdots] $.}
The columns of $\MD=[\Vd_1,\cdots,\Vd_{n_{s}}]$ are the measurements corresponding to source $\Vb_i$.
We have
\begin{eqnarray}\label{minproblem}
  \min\limits_{\Vp} \|\Vr(\Vp)\|_2^2
  & = &
  \min\limits_{\Vp} \sum_{j=1}^{n_{s}} \|\MC^T \MA^{-1}(\Vp) \Vb_j - \Vd_j \|_2^2=
  \min\limits_{\Vp} \|\MC^T \MA^{-1}(\Vp) \MB - \MD\|_F^2. \mbox{  }
\end{eqnarray}
Each evaluation of the objective function requires solving $n_s\cdot n_\omega$ linear systems.
Haber et al.  \cite{HabeChunHerr2012} reduce this cost using simultaneous random sources
combined with (stochastic) trace estimators, following Hutchinson \cite{Hutchinson}.

Let  $\Vw$ be a random vector with mean $\Vo$ and identity covariance matrix,
and let $\mathbb{E}$ denote the expected value with respect to the random vector $\Vw$. Then
\[
  \mathbb{E} \left[ \Vw^T\MRr(\Vp)^T\MRr(\Vp)\Vw \right] =
    \mathrm{trace}  \left( \MRr(\Vp)^T\MRr(\Vp)\right)=\|\MRr(\Vp)\|_F^2.
\]

As a particular choice, we choose $\randsrc$ to be a realization from the Rademacher distribution, where each component of $\randsrc$ is independently and
identically distributed (i.i.d.) taking values from $\{-1, +1\}$, each with probability $\displaystyle \frac{1}{2}$.
Then, as shown in \cite{Hutchinson}, $\randsrc^T\Res(\Vp)^T\Res(\Vp)\randsrc$ is
a minimum variance and unbiased estimator of the trace of $\Res(\Vp)^T\Res(\Vp)$. Thus, the nonlinear least squares problem can be written as a stochastic minimization problem
\begin{eqnarray}\label{stochopt}
 \min\limits_{\Vp} \|\MRr(\Vp)\|_F^2=  \min\limits_{\Vp} \mathrm{trace}  \,\, \MRr(\Vp)^T\MRr(\Vp)=\min\limits_{\Vp} \mathbb{E}  \left( \Vw^T\MRr(\Vp)^T\MRr(\Vp)\Vw \right).
\end{eqnarray}
For a random vector $\Vw$ and simultaneous random source $\MB \Vw$, we have
\begin{equation}\label{Rw}
  \MRr(\Vp)\Vw = (\MC^T \MA^{-1}(\Vp)\MB - \MD)\Vw=\MC^T \MA^{-1}(\Vp)\MB\Vw-\MD\Vw.
\end{equation}
So, computing $||\MRr (\Vp) \Vw||_2^2\ $ requires a single PDE solve rather than $n_s$ solves,
which drastically reduces the cost of a function evaluation.

In contrast to the approach in \cite{HabeChunHerr2012}, we use a Newton-type
method, so we also need to reduce
the cost of Jacobian evaluations.
Therefore, we propose a variation that also drastically reduces the cost of
computing $\MA^{-T}(\Vp)\MC$ for the Jacobian. Let $\Vv \in \mathbb{R}^{n_d}$ and $\Vw \in \mathbb{R}^{n_s}$ with all components
i.i.d.\  uniformly from $\{-1,+1\}$. Using the well-known cyclic 
property of the trace \cite[p. 110]{MeyerMatrixAnalyBook}, we get 
\begin{eqnarray}\label{variation}
  \mathbb{E} \left[ \left(\Vv^T  \MRr \Vw \right)^2\right] & = &
    \mathbb{E}  \left[ \left(\Vv^T \MRr \Vw \right) \left(\Vw^T \MRr^T \Vv \right)\right]
    =\mathbb{E}  \left[ \mathrm{trace} \left(\Vv \Vv^T \MRr \Vw \Vw^T \MRr^T \right)  \right]    \nonumber \\
    &=&
  \mathrm{trace}  \left( \mathbb{E} \left( \Vv \Vv^T \MRr \Vw \Vw^T \MRr^T \right)  \right)
    =\mathrm{trace}  \left( \MRr \MRr^T \right)=||\MRr ||_F^2,
\end{eqnarray}
which requires a single additional adjoint solve rather than an additional $n_d$ solves for the Jacobian. 

Typically, we need multiple random samples $\Vw_j$ and $\Vv_j$ to make the variance
in our stochastic estimates sufficiently small. Hence, we set
\begin{equation}\label{defW}
  \MW = \frac{1}{\sqrt{\ell_s}} (\Vw_1\: \Vw_2\: \cdots \: \Vw_{\ell_{s}}) \;\; \in \Rmn{n_{s}}{\ell_{s}} ,
\end{equation}
where each column vector $\Vw_j$ is i.i.d.\  with zero expectation and covariance equal to the identity
and $\ell_s \ll n_s$.
Similarly, we set
\begin{equation}\label{defV}
  \MV = \frac{1}{\sqrt{\ell_d}}(\Vv_1\: \Vv_2\: \cdots \: \Vv_{\ell_{d}}) \;\;  \in \Rmn{n_{d}}{\ell_{d}} ,
\end{equation}
where each column vector $\Vv_j$ is i.i.d.\  with zero expectation and covariance equal to the identity
and $\ell_d \ll n_d$. It is easily verified that these choices give
\eqs \label{eq:EWWT=I}
  \Exp[\MW \MW^T] = \MI_{n_s} &\mbox{  and  }& \Exp[\MV \MV^T] = \MI_{n_d} .
\eqe
Next, we replace the sources $\MB$ by simultaneous random sources $\MB\Randsrc$ and
the detectors $\MC$ by simultaneous random detectors $\MC\Randdet$. Assume that $\Randsrc$ and $\Randdet$ are independent and we compute unbiased estimates for $\| \Res(\Vp) \|^2$.

\begin{theorem} \label{teo:VRW}
Let $\MW \in \Rmn{n_{s}}{\ell_{s}}$ and $\MV \in \Rmn{n_{d}}{\ell_{d}}$ be as given above.
Let $ \MRr \in \Rmn{n_{d}}{n_{s}}$. Then
\begin{eqnarray}
  \Exp  \left[ \| \MV^T \MRr \MW \|_F^2 \right]= \|\MRr \|_F^2.
\end{eqnarray}
\end{theorem}
\begin{proof}
\begin{eqnarray}
  \mathbb{E} \left[ \| \MV^T \MRr \MW \|_F^2 \right]
  & = &
    \mathbb{E}  \left[ \mathrm{trace} \left(\MW^T \MRr^T \MV \MV^T \MRr \MW \right) \right]
    = \mathbb{E}  \left[ \mathrm{trace} \left(\MW \MW^T \MRr^T \MV \MV^T \MRr \right) \right] \nonumber \\
  & = &
    \mathrm{trace} \left( \mathbb{E} \left[ \MW \MW^T\right] \MRr^T \mathbb{E} \left[\MV \MV^T \right] \MRr \: \right)
      = \mathrm{trace} \left(\MRr^T \MRr \right) = \|\MRr \|_F^2.
\end{eqnarray}
\end{proof}

Since TREGS has proven very effective for the nonlinear least squares problem
in DOT with PaLS, we continue to use the TREGS algorithm in the stochastic minimization
problem
\begin{equation}\label{Eq: true stoch prb}
\min\limits_{\Vp}\mathbb{E}  \left[ || \Randdet^T \Res(\Vp) \Randsrc ||^2 \right] =\min\limits_{\Vp}  ||\Res(\Vp) ||^2.
\end{equation}

We derive the least squares problem used in TREGS to compute a
regularized Gauss-Newton update for the stochastic
problem as follows. For any $\Vp$,
\begin{equation}
\mathrm{vec} \left(\Randdet^T\Res(\Vp)\Randsrc\right) = \left( \MW^T\otimes \MV^T \right) \mathrm{vec}(\MRr(\Vp))
    = \left( \MW^T\otimes \MV^T \right) \Vr (\Vp) ;
\end{equation}
see \cite[lemma 4.3.1]{HornJohnsonMatrixAnaly}.
Using a first order approximation to $\Vr (\Vp+ \Vgd)$ gives
\[
  \left(\MW^T\otimes \MV^T\right)  \Vr (\Vp+\Vgd)\approx \left(\MW^T\otimes \MV^T\right) (\Vr (\Vp)+\MJ \Vgd) ,
\]
which leads to the (sampled) least squares problem
\begin{equation}\label{TREGS J}
  \min\limits_{\pmb{\d}} \|(\MW^T \otimes \MV^T) \MJ {\pmb{\d}}  +(\MW^T \otimes \MV^T) \Vr (\Vp)\|_{2}^2,
\end{equation}
replacing (\ref{Eq:NLLS}). Note that setting up the least squares problem (\ref{TREGS J})
does not require any computations beyond $\MA(\Vp)^{-1} (\MB \MW)$ and $\MA(\Vp)^{-T} (\MC \MV)$.
In addition, (\ref{TREGS J}) has the following desirable properties for the sampled Jacobian and residual,
which follow directly from (\ref{eq:EWWT=I}) and well-known properties of the
Kronecker product.

\eqs
  \label{eq:UnbiasedGrad}
  \Exp\left[ \left( (\MW^T \otimes \MV^T)\MJ \right)^T  \,  ( \MW^T \otimes \MV^T) \Vr \right] & = &
    \MJ^T \Exp[ \MW\MW^T \otimes \MV\MV^T ] \Vr = \MJ^T \Vr  , \\
  \label{eq:UnbiasedHess}
  \Exp\left[ \left( (\MW^T \otimes \MV^T) \MJ \right)^T  \,  (\MW^T \otimes \MV^T)\MJ \right] & = &
    \MJ^T \Exp\left[ \MW\MW^T \otimes \MV\MV^T \right] \MJ = \MJ^T \MJ .
\eqe
So, the proposed randomization provides unbiased estimates for the gradient and the Gauss-Newton Hessian.

Two approaches to stochastic optimization are commonly used \cite{ShapDentRusz2009}. One
approach, stochastic approximation (SA), uses a new random vector (or small batch of random vectors)
in each optimization step. The other approach, sample average approximation (SAA),
uses a fixed set of random vectors over multiple (or many) optimization steps. In this paper, we focus on the SAA approach \cite{ShapDentRusz2009} to solve the stochastic problem
(\ref{Eq: true stoch prb}). The SAA approach approximates (\ref{Eq: true stoch prb}) by the
sample average problem. At each iteration, this approach requires solving only
$\ell_s + \ell_d$ linear systems for each frequency to estimate the objective function and the
Jacobian rather than $n_s + n_d$.

\begin{figure}
    \centering
    \begin{subfigure}[b]{0.32\textwidth}
        \centering
        \includegraphics[width=0.95\textwidth,]{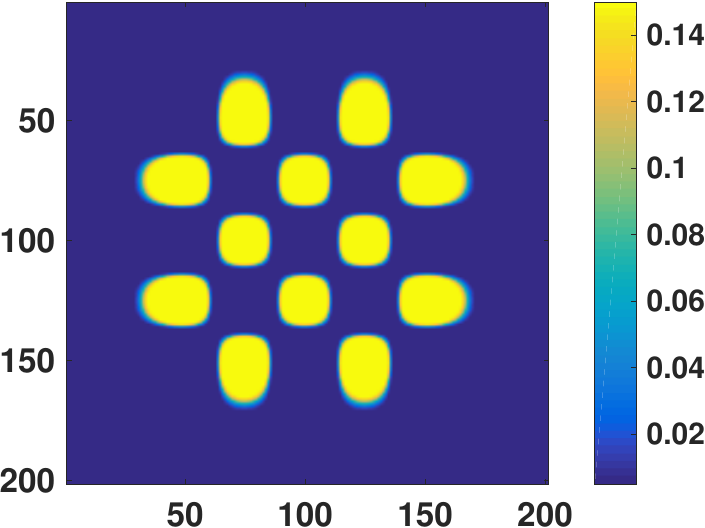}
        \subcaption{}
    \end{subfigure}
        \centering
    \begin{subfigure}[b]{0.32\textwidth}
        \centering
        \includegraphics[width=0.95\textwidth,]{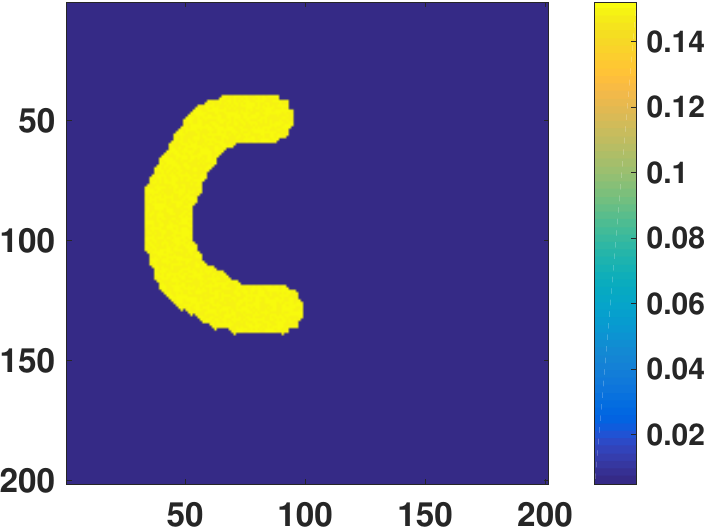}
        \subcaption{}
    \end{subfigure}\\
    \begin{subfigure}[b]{0.32\textwidth}
        \centering
        \includegraphics[width=0.95\textwidth,]{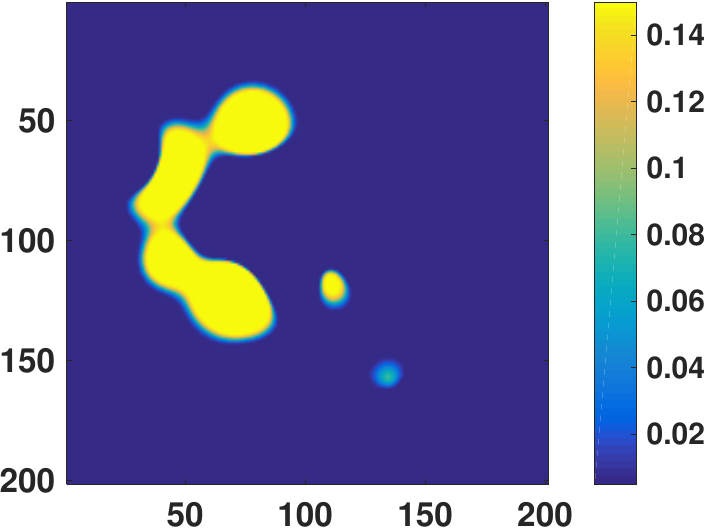}
        \subcaption{}
    \end{subfigure}
    \hfill
    \begin{subfigure}[b]{0.32\textwidth}
        \centering
        \includegraphics[width=0.95\textwidth, ]{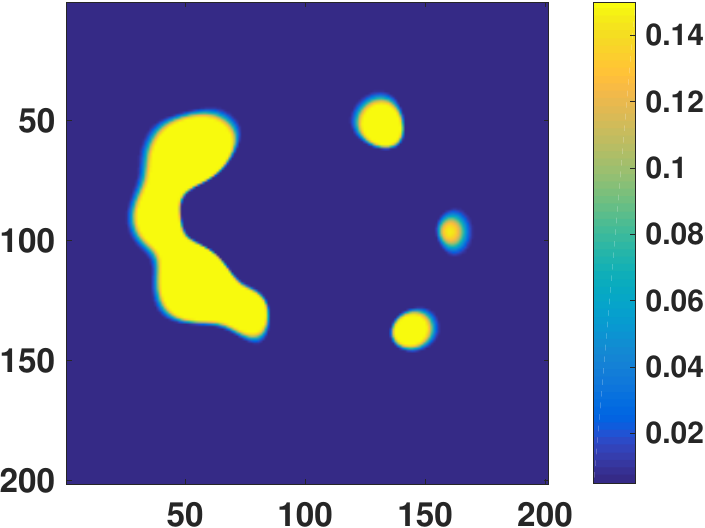}
        \subcaption{}
    \end{subfigure}
        \hfill
      \begin{subfigure}[b]{0.32\textwidth}
        \centering
        \includegraphics[width=0.95\textwidth,]{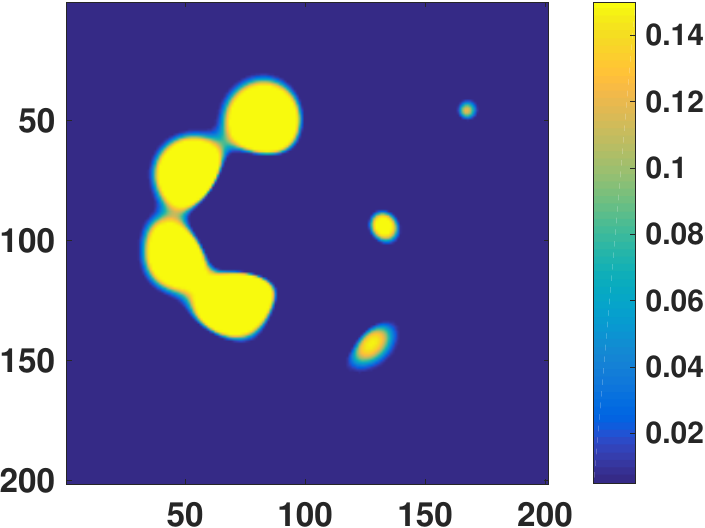}
        \subcaption{}
    \end{subfigure}
\caption{Reconstruction of a test anomaly on $201\times 201$ mesh with 32 sources, 32 detectors, using only the zero frequency.\\
(a) Initial configuration with 25 basis functions arranged in a $5\times 5$ grid where 12 basis functions have positive expansion factors (visible as high absorption regions) and 13 basis functions have negative expansion factors (invisible). (b) True shape of the anomaly. (c) Reconstruction using
all sources and detectors. (d) \& (e) Two reconstruction results using simultaneous random sources and detectors with $\ell_{s}=\ell_{d}=10$.}
\label{Fig:Recons rand}
\end{figure}

We give two representative solutions for our problem using the SAA approach in
Figure \ref{Fig:Recons rand}. For DOT, the use of simultaneous
random sources and detectors initially leads to good progress.
However, later in the iteration the convergence slows down,
and in many cases, for our problem, it does not lead to sufficiently accurate solutions.
In fact, with the SAA approach, (typically) the residual norm
does not reach the noise level, the stopping criterion used
(see the discussion of PaLS in
Section~\ref{Sec:SimultOpt_backg}), while
the standard optimization using all sources and all detectors
does converge to the noise level.
We will demonstrate this in Section \ref{Sec:SimultOpt_NumExp}.
In the next section, we provide a solution to this problem.

\subsection{Improving the Randomized Approach}\label{sec:OptSrcDet_new}
In the standard SAA approach, when convergence slows down
or a minimum is found for the chosen sample (but not for the
true problem), a new sample is chosen to improve the approximate solution.
However, for our problem this approach leads to slow convergence and stagnation,
unless we obtain more accurate estimates of the dominant singular components of the
Jacobian and the corresponding components of the gradient (see also \cite{StuKil11c}).
Hence, after exploiting the relatively
fast initial convergence for our problem, we want to avoid stagnation of convergence in the next phase.
One approach is to add additional random simultaneous sources
and detectors, that is, increase the sample size, as proposed in \cite{FarbodDoelAscher2,ByrdChinNocedalWu, BollaByrdNocedal, FarbodMahoney1, FarbodMahoney2}, with good results.
However, this requires progressively more, expensive, solves.
Therefore, for efficiency, we choose to keep the number of simultaneous sources
and detectors fixed.
To improve convergence, we exploit the (often fairly good)
approximate solution at a chosen modest intermediate tolerance
to obtain a small number of simultaneous sources
and detectors that are optimized to provide more accurate
estimates of the desired derivative information.
Since this optimization is local for the current $\Vp$,
we complement these optimized vectors by random simultaneous
sources and detectors. We make this precise below.

The nonlinear least squares algorithm TREGS focuses on the dominant singular
values of the Jacobian to compute good updates to the parameter vector \cite{StuKil11c}.
The corresponding right singular vectors capture the directions in parameter
space of largest sensitivity in the objective function.
Hence, we want to
update $\MW$ and $\MV$ so as to
best approximate the dominant right singular subspace
of $\MJ$ while respecting the Kronecker product structure
in (\ref{TREGS J}).
This is important for two reasons. First, for the same (fixed) small number of simultaneous sources and detectors, this gives us
locally (at the current $\Vp$) the best approximation to what TREGS would do using all sources
and detectors. Second, the directions corresponding to the dominant right singular vectors
are best informed by the data.

So, when a chosen intermediate
tolerance is reached, our method computes
the full Jacobian $\MJ$ once, which requires a total of $n_s+n_d$ solves. Then,
we compute a small number, $q_s$ respectively $q_d$, of
orthonormal, optimized simultaneous sources ($\wh{\Randsrc}$)
and detectors ($\wh{\Randdet}$). In practice, small $q_\mathrm{s}$ and $q_\mathrm{d}$, 2 to 4, seem to be
sufficient. We provide some experimental
results regarding the number of optimized directions in Section~\ref{Sec:SimultOpt_NumExp}.
Since $\MJ$ is typically of rank substantially lower than $n_\mathrm{p}$
\cite{Aghasi_etal11,StuKil11c}, we expect that computing optimized directions
can be done with a cheap approximation to $\MJ$, for example,
computed using the  that can be computed
at substantially lower cost than solving for all sources and detectors (for each
frequency) \cite{HalMartTropp}.
However, this is beyond the scope of the current paper.

We would like to maximize
\eqs\label{eq:maxprob}
  \|(\wh{\MW}^T \otimes \wh{\MV}^T) \MJ \|_F^2 .
\eqe
However,
the Kronecker product structure combined with the constraints that
$\wh{\MW}$ and $\wh{\MV}$ be isometric matrices leads to a nonlinear constrained
optimization problem. For efficiency, we exploit the tensor-structure
of this problem, that is, we consider the Jacobian as a third-order array,
$\cJ$, with components
\eqs \label{eq:Jtensor}
  \cJ_{ijk} & = &
  \Vc_i^T \MA^{-1}\frac{\partial \MA}{\partial p_k} \MA^{-1} \Vb_j .
\eqe
This allows us to use an alternating least squares algorithm,
a variant of the Higher Order Orthogonal Iteration (HOOI) \cite{LatMooVan00b}, to find
$\wh{\MV} \in \Rmn{n_\mathrm{d}}{q_\mathrm{d}}$ and
$\wh{\MW} \in \Rmn{n_\mathrm{s}}{q_\mathrm{s}}$ that
approximately maximize (\ref{eq:maxprob}); details follow in the next section.
The algorithm is only guaranteed to find a local maximum \cite{KroLee80,LatMooVan00b}.
However, in a number of numerical tests carried out, the tensor algorithm discussed below
seems to always converge to the global maximum. A similar observation is
reported in \cite{LatMooVan00b}. In our experiments, the algorithm also attains
the same solutions as \textsc{MATLAB}\textsuperscript{\textregistered}'s {\tt fmincon} routine, which optimizes for $\wh{\MW}$ and $\wh{\MV}$ simultaneously.

In the remainder of this section,
we first discuss computing the optimized simultaneous
detectors ($\wh{\MV}$) and sources ($\wh{\MW}$) and then
complementing these with random simultaneous sources
constrained to $\Ra{\wh{\MW}}^{\perp}$ and
random simultaneous detectors
constrained to $\Ra{\wh{\MV}}^{\perp}$.

\subsubsection{Computing Optimized Simultaneous Sources and Detectors}
This problem is closely related to the
truncated higher-order SVD (HOSVD) \cite{LatMooVan00a,KoldBad09} or,
more precisely, to a truncated Tucker2 decomposition \cite{Tuck66,KoldBad09},
as we do not need truncation in the parameter-derivative direction (the
columns of $\MJ$).
As $\wh{\MW}$ and $\wh{\MV}$ are both isometric matrices, so is their
Kronecker product, $\MX = \wh{\MW} \otimes \wh{\MV}$. Let
$[ \MX \; \MX_c ]$ be a (real) orthogonal matrix, and let
$S^{ k \times \ell} = \{ \MgVC \in \Rmn{k}{\ell} | \: \MgVC^T\MgVC = \MI_{\ell} \}$
(the set of all $k \times \ell$ isometric matrices).
Then, it follows from standard
properties of the Frobenius norm that
$\| \MJ \|_F = \| [\MX \; \MX_c]^T \MJ \|_F$ and
\eqsn
  \| \MJ \|_F^2 & = & \| \MX^T \MJ \|_F^2 + \| \MX_c^T \MJ \|_F^2 .
\eqen
Since $\| (\MI - \MX \MX^T) \MJ \|_F = \| \MX_c^T \MJ \|_F$,
we have that the maximization problem
\eqs\label{eq:BestApproxEquiv}
  \arg \max_{
    \oars{c}
      \wt{\MW} \in \MS^{n_s \times \ell_s} \\
      \wt{\MV} \in \MS^{n_d \times \ell_d}
    \oare
    }
  \|(\wt{\MW}^T \otimes \wt{\MV}^T) \MJ \|_F^2
\eqe
is equivalent with the minimization problem
\eqsn
  \arg \min_{
    \oars{c}
      \wt{\MW} \in \MS^{n_s \times \ell_s} \\
      \wt{\MV} \in \MS^{n_d \times \ell_d}
    \oare
    }
  \| \MJ - (\wt{\MW} \otimes \wt{\MV}) (\wt{\MW} \otimes \wt{\MV})^T \MJ \|_F^2 ;
\eqen
see also \cite[Theorem 4.1]{LatMooVan00b} and \cite[p. 477]{KoldBad09}.
We solve the maximization problem by a slight variation of the
HOOI algorithm \cite{LatMooVan00b} to
compute a truncated Tucker2 decomposition
\cite{Tuck66,KoldBad09}. We do not need any compression of the
column dimension of $\MJ$ (the direction
corresponding to the parameter derivatives).

We briefly outline the algorithm in tensor form; a
\textsc{MATLAB}\textsuperscript{\textregistered}
pseudo-code is
given in subsection~\ref{subsec:Implementation}.
Let $\cJ$ be the tensor representation of the Jacobian $\MJ$, where
$\cJ_{ijk} = \Vc_i^T \MA^{-1}\frac{\partial \MA}{\partial p_k} \MA^{-1} \Vb_j$;
the dependence on $\Vp$ in $\MA$, $\MJ$, and $\cJ$ is suppressed for
brevity.

\el
(1) We compute the SVD of the matrix obtained by unfolding
the tensor $\cJ$ along its lateral slices $\cJ_{:,j,:}$ and computing the
SVD of the resulting matrix $\wt{\MJ}$.
\eqs \label{eq:Jtilde}
  \MU \MgSi \MY^T = \wt{\MJ} & = &
  \ars{cccc}
    \cJ_{:,1,:} & \cJ_{:,2,:} & \cdots & \cJ_{:,n_s,:}
  \are .
\eqe
Next, we set the initial $\wt{\MV} = [ \Vu_1 \; \Vu_2 \; \ldots \; \Vu_{q_d}]$
(the leading $q_d$ left singular vectors of $\wt{\MJ}$).

Subsequently, we iterate the following two steps until the
change in the approximate solution, the tuple $(\wt{\MV}, \wt{\MW})$, is sufficiently small.
\el
(2) Let the tensor $\cK$ be defined as
\[
  \cK_{jk \ell} = \sum_{m = 1}^{n_d}
    (\wt{\Vv}_{\ell})_m \cJ_{m,j,k} ,
    \qquad \mbox{for } \ell = 1, 2, \ldots, q_d ,
\]
where the vectors $\wt{\Vv}_{\ell}$ are the columns of $\wt{\MV}$
($\ell = 1, \ldots, q_d$).
We compute the SVD of the matrix $\wt{\MK}$ obtained by unfolding $\cK$
along its frontal slices $\cK_{:,:,\ell}$,
\eqs \label{eq:Ktilde}
  \MgF \MgO \MgY^T = \wt{\MK} & = &
  \ars{cccc}
    \cK_{:,:,1} & \cK_{:,:,2} & \cdots & \cK_{:,:,q_d}
  \are ,
\eqe
and we set $\wt{\MW} = [ \Vgf_1 \; \Vgf_2 \; \ldots \; \Vgf_{q_s}]$.
(the leading $q_s$ left singular vectors of $\wt{\MK}$)

\el

(3) Define the tensor
$\cL$ by
\[
  \cL_{ik \ell} = \sum_{m=1}^{n_s} (\wt{\Vw}_{\ell})_m \cJ_{i,m,k} ,
  \qquad \mbox{for } \ell = 1, 2, \ldots, q_s ,
\]
where the vectors $\wt{\Vw}_{\ell}$ are the columns of
$\wt{\MW}$ ($\ell = 1, \ldots, q_s$).
We compute the SVD of the matrix $\wt{\ML}$ obtained by unfolding
$\cL$ along its frontal slices $\cL_{:,:,\ell}$,
\eqs \label{eq:Ltilde}
  \MU \MgSi \MY^T = \wt{\ML} & = &
  \ars{cccc}
    \cL_{:,:,1} & \cL_{:,:,2} & \cdots & \cL_{:,:,q_s}
  \are .
\eqe
and we set (the new) $\wt{\MV} = [ \Vu_1 \; \Vu_2 \; \ldots \; \Vu_{q_d} ]$.

\el
Finally, after convergence, we set $\wh{\MV} = \wt{\MV}$ and $\wh{\MW} = \wt{\MW}$.

While this procedure gives good approximate solutions,
in general only a local maximum of (\ref{eq:maxprob}) is guaranteed
\cite{LatMooVan00b}.
However, in a number of numerical test carried out, the tensor algorithm
seems to always converge to the globally optimal solution.
This experience was also reported in \cite{LatMooVan00b}.
Moreover, we have the following useful result for the case that $\| \MJ \|_F$ can be preserved
exactly.
\teos
If isometric matrices $\wt{\MW} \in \Rmn{n_s}{q_s}$ and $\wt{\MV} \in \Rmn{n_d}{q_d}$ exist
such that
\eqs \label{eq:opt_sol}
  \| (\wt{\MW}^T \otimes \wt{\MV}^T) \MJ \|_F & = & \| \MJ \|_F,
\eqe
then, in one iteration, steps (1) -- (3) above compute
isometric matrices
$\wh{\MW} \in \Rmn{n_s}{q_s}$ and $\wh{\MV} \in \Rmn{n_d}{q_d}$ such that
\eqsn
  \| (\wh{\MW}^T \otimes \wh{\MV}^T) \MJ \|_F & = & \| \MJ \|_F .
\eqen
\teoe
\prfs
This result follows from the error bound (Property 10) in \cite{LatMooVan00a},
which in the (matrix) notation of this paper is given by
\eqs \label{eq:TensorApproxBnd}
  \|\MJ - (\wt{\MW} \otimes \wt{\MV})(\wt{\MW} \otimes \wt{\MV})^T \MJ\|_F^2
  & \leq &
  \sum_{i = q_d+1}^{n_d} \s_i^2 + \sum_{j = q_s+1}^{n_s} \o_j^2 .
\eqe
The assumption (\ref{eq:opt_sol}) implies that
$\Ra{\wt{\MJ}} \subseteq \Ra{\wt{\MV}}$, which in turn implies that
$\Rk{\wt{\MJ}} \leq q_d$ and hence that
$\s_{q_d+1} = 0, \s_{q_d+2} = 0, \ldots, \s_{n_d}=0$.
In a similar fashion, from assumption (\ref{eq:opt_sol}) combined with the
choice for $\wt{\MV}$ from step (1), see (\ref{eq:Jtilde}), it follows that
$\Ra{\wt{\MK}} \subseteq \Ra{\wt{\MW}}$, which implies that
$\Rk{\wt{\MK}} \leq q_s$ and hence
$\o_{q_s+1} = 0, \o_{q_s+2} = 0, \ldots, \o_{n_s}=0$.
\newline\indent
Substitution of $\s_{q_d+1}, \ldots, \s_{n_d}, \o_{q_s+1}, \ldots, \o_{n_s}$ into
(\ref{eq:TensorApproxBnd}) shows that the algorithm will have converged
at this point. The algorithm then sets $\wh{\MW} = \wt{\MW}$ and
$\wh{\MV} = \wt{\MV}$.
\prfe

Since $\wh{\MW}$ and $\wh{\MV}$ are only
optimal at the current $\Vp$, we complement these optimized simultaneous sources
and detectors with a new set of random simultaneous sources and detectors constrained to
the orthogonal complement of the span of the optimized directions,
keeping the total number of columns in $\MW$ and $\MV$ the same as before.
This procedure can be carried out periodically or for a sequence of prescribed tolerances,
but in our experiments it never needs to be done more than once.

\subsubsection{Computing Complementary Random Simultaneous Sources and Detectors}
We extend the optimized sources and detectors with random simultaneous
sources and detectors.
Let $\MW_f = [ \wh{\MW} \, \MW_c] \in \Rmn{n_\mathrm{s}}{n_\mathrm{s}}$ and
$\MV_f = [ \wh{\MV} \, \MV_c] \in \Rmn{n_\mathrm{d}}{n_\mathrm{d}}$ be orthogonal
matrices. We have
\eqs
  \nonumber
  \| \MRr(\Vp) \|_F & = & \| \MV_f^T \MRr(\Vp) \MW_f \|_F  \\
  & = &
  \left\|
  \ars{cc}
    \wh{\MV}^T \MRr(\Vp) \wh{\MW}  &  \wh{\MV}^T \MRr(\Vp) \MW_c     \\
    \MV_c^T \MRr(\Vp) \wh{\MW}        &  \MV_c^T \MRr(\Vp) \MW_c
  \are
  \right\|_F
  \label{eq:VfRWf_1}
\eqe
The $(1,1)$-block of this matrix can be computed using the known optimized
sources and detectors. We estimate the remaining blocks, proceeding more or
less as before.
We pick random matrices
$\MY = (\ell_s - q_s)^{-1/2}[\Vy_1 \, \Vy_2 \, \ldots \, \Vy_{\ell_s - q_s}]$, where each
column vector $\Vy_j \in \Rn{{\ell_s - q_s}}$ is i.i.d.\  with zero mean and identity covariance and
$\MZ = (\ell_d - q_d)^{-1/2}[ \Vz_1 \, \Vz_2 \, \ldots \, \Vz_{\ell_d - q_d} ] $, where each
column vector
$\Vz_j \in \Rn{{n_d - q_d}}$ is i.i.d.\  with zero mean and identity covariance.
In our numerical experiments, we again use the Rademacher distribution. 
Next, we set the new matrices $\MW$ and $\MV$ to
\eqs\label{eq:newW}
  \MW & = & [ \wh{\MW} \; (\MW_c\!\MY) ] , \\
  \label{eq:newV}
  \MV & = & [ \wh{\MV} \; (\MV_c\,\MZ) ] .
\eqe
We have the following results.
\teos \label{teo:UnbiasComplRnd}Let $\Randsrc \in \Rmn{n_{s}}{\ell_{s}}$ and $\Randdet \in \Rmn{n_{d}}{\ell_{d}}$ be given in (\ref{eq:newW}) and (\ref{eq:newV}), respectively.
Let $ \Res \in \Rmn{n_{d}}{n_{s}}$. Then,
\eqs
  \label{eq:covW}
  \Exp[ \MW \MW^T ] & = &  \MI_{n_s}  ,   \\
  \label{eq:covV}
  \Exp[ \MV \MV^T  ]   & = & \MI_{n_d}  ,  \\
  \label{eq:newVRW}
  \Exp[ \| \MV^T \MRr(\Vp) \MW \|_F^2 ] & = & \| \MRr(\Vp) \|_F^2 .
\eqe
\teoe
\prfs
For (\ref{eq:covW}), we have
\eqsn
  \Exp \left[ [ \wh{\MW} \; (\MW_c\!\MY) ] [\wh{\MW} \; (\MW_c\!\MY) ]^T \right]  & = &
    \wh{\MW}\wh{\MW}^T + \Exp[ \MW_c \MY \MY^T \MW_c^T ] \\
  & = &
  \wh{\MW}\wh{\MW}^T +  \MW_c \Exp [ \MY \MY^T ] \MW_c^T = \MI_{n_s} .
\eqen
An analogous derivation holds for (\ref{eq:covV}).
The proof for the last equation follows the proof for
Theorem~\ref{teo:VRW}, using the results above.
\prfe

As a result of Theorem~\ref{teo:UnbiasComplRnd}, the new $\MW$ and $\MV$
again give for the expectation of the sampled gradient and the expectation of the
sampled Gauss-Newton
Hessian the true gradient and Gauss-Newton Hessian,
\eqsn
  \Exp[ ( (\MW^T \otimes \MV)\MJ )^T (\MW^T \otimes \MV) \Vr ] & = &
     \MJ^T \Vr  ,  \\
  \Exp[ ( (\MW^T \otimes \MV^T)\MJ )^T( (\MW^T \otimes \MV^T)\MJ) ] & = &
     \MJ^T \MJ .
\eqen

\subsection{Implementation}\label{subsec:Implementation}
In this section, we first give an algorithm for the computation of optimized simultaneous sources and detectors. Algorithm~\ref{alg:opt_srcdet} is based on the matrix representation of the Jacobian as defined in (\ref{eq:Jacobiancomp}). Next, we outline the efficient computation of the residual and Jacobian.
\begin{algorithm}[t]
\caption{Compute Optimized Sources and Detectors}
\label{alg:opt_srcdet}
\begin{algorithmic}
\State {\em \{ compute $\wt{\MJ}$ and its $q_d$ leading left singular vectors \}}
\State $\wt{\MJ}$ = [\:]; \Comment {initialize $\wt{\MJ}$ as empty matrix}
\For{ j = 1:$n_s$}
  \State $\wt{\MJ} = [\wt{\MJ} \; \MJ_{(j-1)n_d+1:jn_d,1:k}]$
  \Comment{add next block from tensor $\cJ$}
\EndFor
\State [$\MU$,$\MgSi$,$\MY$] = SVD($\wt{\MJ}$);\:
  $\wt{\MV}$ = [$\Vu_1$ $\Vu_2$ ... $\Vu_{q_d}$]
\Comment{keep $q_d$ leading left singular vectors of $\wt{\MJ}$}
\While{not converged}
  \State {\em \{ compute $\wt{\MK}$ and its $q_s$ leading left singular vectors \} }
  \State $\wt{\MK}$ = [\:]
  \Comment{initialize $\wt{\MK}$ as empty matrix}
  \For{$\ell$ = 1:$q_d$}
    \State $\wt{\MK}_{\ell}$ = [\:]
    \For{j = 1:$n_s$}
      \State $(\wt{\MK}_{\ell})_{j,1:k} = \Vv_{\ell}^T \MJ_{(j-1)n_d+1:jn_d,1:k}$
      \Comment{compute row $j$ of $\wt{\MK}_{\ell}$}
    \EndFor
    \State $\wt{\MK} = [\wt{\MK} \; \wt{\MK}_{\ell}]$
    \Comment{add next block of $\wt{\MK}$}
  \EndFor
  \State [$\MgF$,$\MgO$,$\MgY$] = SVD($\wt{\MK}$);
    $\wt{\MW} = [\Vgf_1 \,\Vgf_2 \, \ldots \, \Vgf_{q_s}]$
  \Comment{keep $q_s$ leading left singular vectors of $\wt{\MK}$}
  \State {\em \{ compute $\wt{\ML}$ and its $q_s$ leading left singular vectors \} }
  \State $\wt{\ML}$ = [\:]
  \For{$\ell$ = 1:$q_s$}
    \State $\wt{\ML}_{\ell} = \sum_{j=1}^{n_s} (\wt{\Vw}_{\ell})_j
      \MJ_{(j-1)n_d+1:jn_d,1:k}$
    \State $\wt{\ML} = [\wt{\ML} \; \wt{\ML}_{\ell}]$
    \Comment{add next block of $\wt{\ML}$}
  \EndFor
  \State [$\MU$,$\MgSi$,$\MY$] = SVD($\wt{\ML}$);
    $\wt{\MV}$ = [$\Vu_1$ $\Vu_2$ ... $\Vu_{q_d}$]
  \Comment{keep $q_d$ leading left singular vectors of $\wt{\MJ}$}
\EndWhile
\State $\wh{\MV} = \wt{\MV}; \; \wh{\MW} = \wt{\MW}$
\end{algorithmic}
\end{algorithm}
We estimate the norm of the residual using
\begin{equation}\label{EsTResidual}
 (\MW^T \otimes \MV^T) \Vr(\Vp)=   \begin{bmatrix}
  \MV^T\MC^T \MA^{-1}(\Vp)\MB\Vw_1 - \MD\Vw_1\\
          \vdots \\
  \MV^T\MC^T \MA^{-1}(\Vp)\MB\Vw_{\ell_{s}}  - \MD\Vw_{\ell_{s}}
  \end{bmatrix}=  \begin{bmatrix}
  \MV^T\MC^T \Vz_1 - \MD\Vw_1\\
          \vdots \\
  \MV^T\MC^T \Vz_{\ell_{s}}  - \MD\Vw_{\ell_{s}}
  \end{bmatrix},
\end{equation}
where we solve $ \MA(\Vp)\Vz_i=\MB\Vw_i$ for $\Vz_i,\: i=1\cdots \ell_s$. This reduces the number of large solves from $n_s$ to $\ell_s$ per frequency. To compute the Jacobian, we solve the systems, $ \MA^T(\Vp)\Vy_j=\MC\Vv_j$ for $\Vy_j,\:j=1\cdots \ell_d$. This reduces the additional number of large solves from $n_d$ to $\ell_d$ per frequency. We can use iterative solvers or sparse direct solvers depending on the size of the system \cite{KilStu2006}. To obtain the $k$-th column of $(\MW^T \otimes \MV^T) \MJ$, we compute
\eqs \label{Jacobian kth column}
 \ars{l}
   \Vy_1^T  \frac{\partial \MA(\Vp)}{\partial \Vp_k}\Vz_1 \cdots \Vy_{\ell_{d}}^T  \frac{\partial \MA(\Vp)}{\partial \Vp_k}\Vz_1 \:\: \Vy_1^T  \frac{\partial \MA(\Vp)}{\partial \Vp_k}\Vz_2  \cdots \Vy_{\ell_{d}}^T  \frac{\partial \MA(\Vp)}{\partial \Vp_k}\Vz_{\ell_{s}}\\
  \are^T,
\eqe
where
$\partial \MA(\Vp) / \partial \Vp_k$
is a diagonal matrix if we only invert for absorption. If we also invert for diffusion, this matrix has 5 (in 2D) or 7 (in 3D) diagonals. Moreover, after  a few iterations, the changes in $\MA(\Vp)$ are highly localized, and
$\partial \MA(\Vp) / \partial \Vp_k$ contains mostly zero coefficients; see \cite{StuGuKilChatBeatOCon}. In that case, we first find the few nonzero components of 
$\partial \MA(\Vp) / \partial \Vp_k$
for each $k$ and the corresponding nonzeros in $\Vz_i$ and $\Vy_j$, and
next we efficiently compute
$\Vy_i^T  \frac{\partial \MA(\Vp)}{\partial \Vp_k}\Vz_j$
referencing only the few nonzero components in
$(\partial \MA(\Vp) / \partial \Vp_k) \Vz_j$.

\section{Numerical Experiments}\label{Sec:SimultOpt_NumExp}
In this section, we provide two numerical experiments, a 2D and a 3D test case, to demonstrate the effectiveness of combining random simultaneous sources and detectors with optimized simultaneous sources and detectors. 
The results show that our approach not only produces reconstruction results 
that are close to those obtained using all sources and all detectors, 
but it also substantially reduces the computational cost.

\begin{figure}
        \centering
       \begin{minipage}{0.32\textwidth}
        \centering
  \includegraphics[width=0.95\textwidth,]{Initial}
        \subcaption{}
    \end{minipage}\:\:\:\:\:\:\:\:\:\:\:\:\:\
      \begin{minipage}{0.32\textwidth}
        \centering
  \includegraphics[width=0.95\textwidth,]{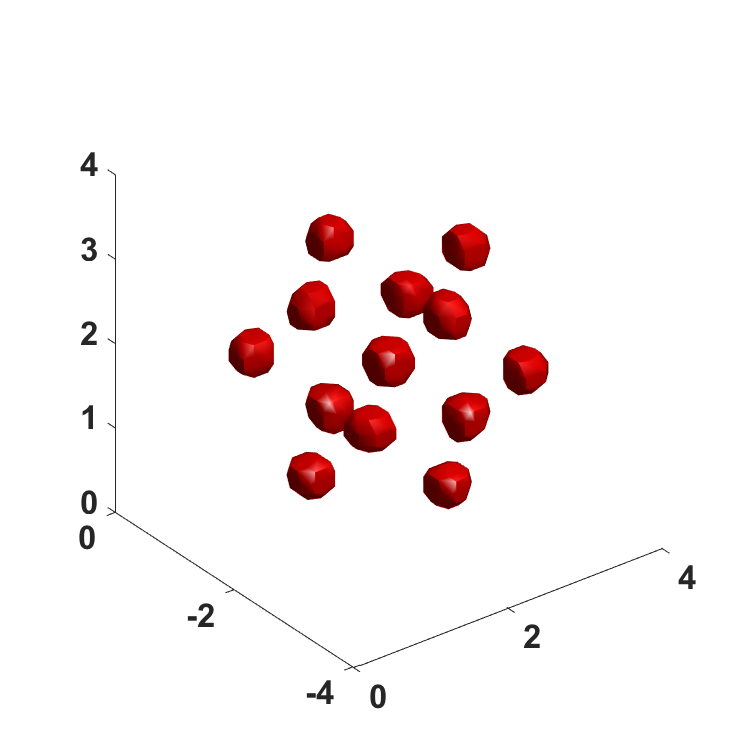}
        \subcaption{}
    \end{minipage}
    \caption{(a) Initial configuration for the 2D experiment with 
    25 basis functions arranged in a $5\times 5$ grid, where 12 basis functions have positive expansion factors (visible) and 13 basis functions have negative expansion factors (invisible). (b) Initial configuration for 
    the 3D experiment with 27 basis functions arranged in a $3\times 3 \times 3$ grid, 
    where 13 basis functions have positive expansion factors (visible) and 14 basis functions have negative expansion factors (invisible).}
 \label{Fig: Initial Conf}
\end{figure}

Our experimental set up is the same as that in \cite{StuGuKilChatBeatOCon}, 
where model reduction was used to reduce the cost of inversion in DOT. 
The absorption images for the initial sets of parameters for the 2D and 3D experiments
are given in Figure~\ref{Fig: Initial Conf}. 
For each test case, we construct anomalies in the pixel basis, and 
we add a small normally distributed random heterogeneity to both the background and 
the anomaly to make the medium inhomogeneous. 
We use this pixel-based absorption image to compute the 
(true) measured data and add $\delta=0.1\%$ white noise to the measured data. 
This is the same noise level as used in \cite{StuGuKilChatBeatOCon}. 
Then, we reconstruct the absorption images using the PaLS \cite{Aghasi_etal11} 
repesentation and TREGS \cite{StuKil11c} for optimization.\\

\textbf{2D Experiment.} We use a 201$\times$201 grid, which yields $40,401$ unknowns
in the discretized PDE (\ref{eq: PDE}). The model has 32 sources, 32 detectors, and we use only the zero frequency. Our model has 25 CSRBFs, which leads to 100 parameters (four per 2D basis function) for the nonlinear optimization. We use the same starting guess for each trial (see Figure~\ref{Fig: Initial Conf}a),
with the 25 basis functions arranged in a $5\times 5$ grid, where 12 basis functions have a positive expansion coefficient (visible as high absorption regions) and 13 basis functions have a negative expansion coefficient (invisible).

We use 10 random simultaneous sources and detectors. We update 
the simultaneous sources and detectors as discussed in Section \ref{Sec:SAA_Opt} 
after a chosen intermediate tolerance has been reached. We find that, in general, the 
the noise level, $\delta$, is a good choice as the intermediate tolerance:  $||\res(\Vp)||_2^2=\delta$. 
Since the PaLS representation regularizes the problem,
we consider the problem converged when $||\res(\Vp)||_2^2 \leq \delta^2$. 
This is called the discrepancy principle (the factor $\frac{1}{2}$ in (\ref{eq:NonlinLS}) is dropped for convenience). We run the 2D experiment for 50 trials. In each trial, the random simultaneous sources and detectors are chosen independently to get representative reconstruction results.

\begin{figure}
    \centering
    \begin{subfigure}[b]{0.3\textwidth}
        \centering
          \includegraphics[width=0.95\textwidth,]{Orig_cupNoFrq}
        \subcaption{}
    \end{subfigure}
        \centering
    \begin{subfigure}[b]{0.3\textwidth}
        \centering
        \includegraphics[width=0.95\textwidth]{TREGS_cupNoFrq}
        \subcaption{}
    \end{subfigure}
    \begin{subfigure}[b]{0.3\textwidth}
        \centering
         \includegraphics[width=0.95\textwidth]{disb_cup_401C_1Sample__Stoch}
        \subcaption{}
    \end{subfigure}
    \centering
       \begin{subfigure}{0.3\textwidth}
        \centering
        \includegraphics[width=0.95\textwidth]{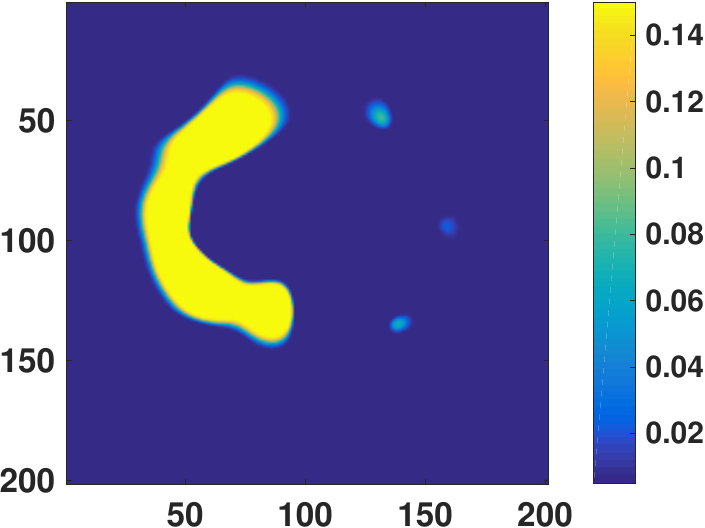}
        \subcaption{}
    \end{subfigure}
    \begin{subfigure}{0.3\textwidth}
        \centering
        \includegraphics[width=0.95\textwidth]{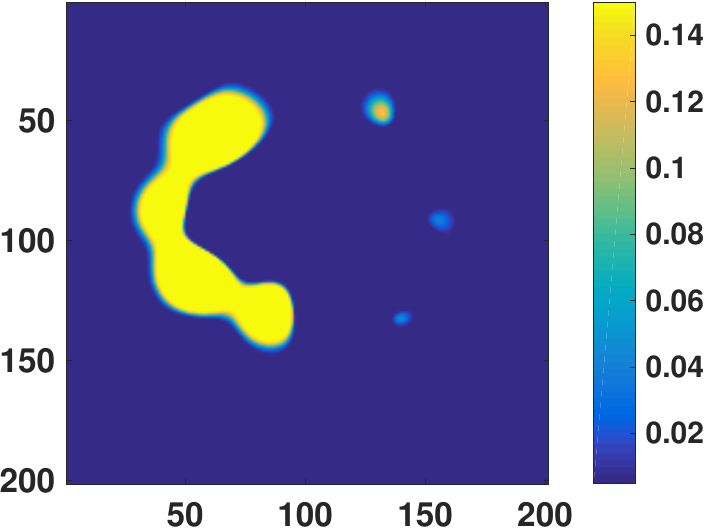}
        \subcaption{}
    \end{subfigure}
      \begin{subfigure}{0.3\textwidth}
        \centering
        \includegraphics[width=0.95\textwidth]{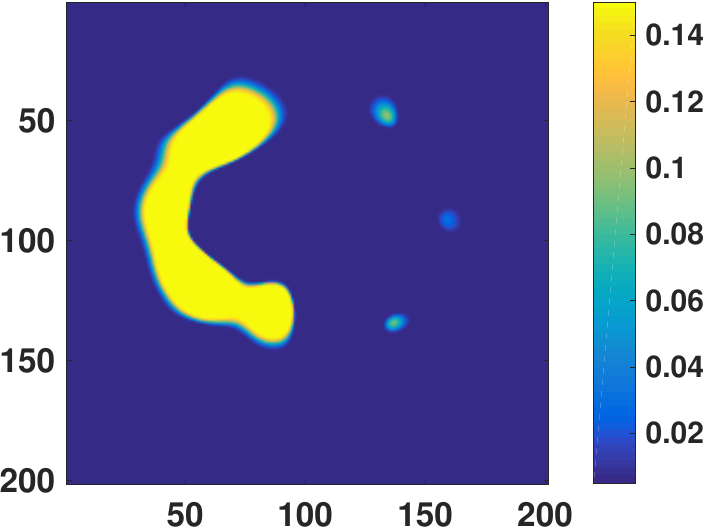}
        \subcaption{}
    \end{subfigure}
\caption{Results for Example 1. Reconstruction of a test anomaly on $201\times 201$ mesh with 32 sources and detectors, 25 basis functions, and using only the zero frequency. The SAA approach uses 10 random simultaneous sources and detectors.\\
(a) True shape of the anomaly. (b) Reconstruction using all sources and all detectors. (c) Reconstruction using the SAA approach at a chosen intermediate tolerance.
(d) Reconstruction with SAA and 1 optimized simultaneous source and detector. (e) Reconstruction with SAA and 2 optimized simultaneous sources and detectors. (f) Reconstruction with SAA and 3 optimized simultaneous sources and detectors.}
\label{Fig:Recons All src det vs SAA vs Repl Exp1_newOpt}
\end{figure}

\begin{figure}
        \centering
       \begin{minipage}{0.32\textwidth}
        \centering
        \includegraphics[width=0.95\textwidth]{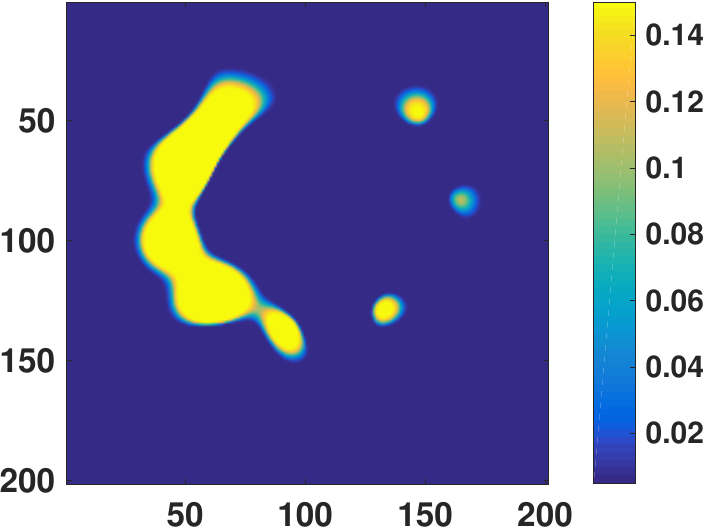}
        \subcaption{}
    \end{minipage}
    \centering
      \begin{minipage}{0.32\textwidth}
        \centering
        \includegraphics[width=0.95\textwidth]{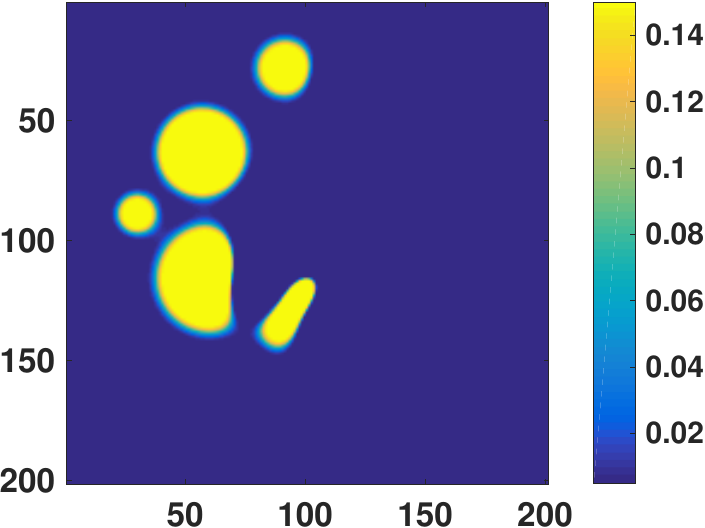} 
        \subcaption{}
    \end{minipage}
\caption{Two examples of poor SAA reconstructions for the 2D problem after the maximum number of iterations. Reconstruction of a test anomaly on $201\times 201$ mesh with 32 sources and detectors, 25 basis functions, and using only the zero frequency.}
\label{Fig: Poor SAA for all exp_newOpt}
\end{figure}

\textbf{Example 1.} The true absorption image for Example 1 is given in Figure~\ref{Fig:Recons All src det vs SAA vs Repl Exp1_newOpt}a. We also include the reconstruction results using all sources and detectors for comparison (see Figure~\ref{Fig:Recons All src det vs SAA vs Repl Exp1_newOpt}b). As can be seen in Figure~\ref{Fig:Recons All src det vs SAA vs Repl Exp1_newOpt}c at the intermediate tolerance, SAA gives a good localization of the anomaly; however, there is little further improvement using SAA (see Figure~\ref{Fig:Recons rand}d-e and Figure~\ref{Fig: Poor SAA for all exp_newOpt}a). Figure~\ref{Fig:Recons All src det vs SAA vs Repl Exp1_newOpt}d-f show that using optimized simultaneous sources and detectors leads to solutions of the same quality as obtained using all sources and detectors. We report the total number of PDE solves required for each approach in Table~\ref{Table: All the numbers for Exp1_newOpt} for a representative result from 50 trials.

While initially the SAA estimate is unbiased,
a systematic underestimation of the residual/misfit (bias) \cite[Section 5.1.2]{ShapDentRusz2009}
arises, since we optimize for a specific small set of random simultaneous sources and detectors. 
As a result, the algorithm generally stops prematurely. This can make a big difference, since 
often substantial improvement in the shape of the anomaly occurs towards
the end of the optimization. Figure \ref{Fig: Poor SAA for all exp_newOpt} demonstrates how poor
the reconstructions using only the SAA approach can be at the convergence tolerance when
underestimation of the residual norm is severe. 
To make a fair comparison in terms of the number of large systems
solved, we check the true function evaluation of the SAA approach on the side.
Table~\ref{Table: Showing biased/unbiased_newOpt} shows that in terms of the true function
evaluation, the SAA approach does not reach the convergence tolerance. 
Once we use a few optimized sources and detectors, this is no longer an issue (see Table~\ref{Table: Showing biased/unbiased_newOpt}).

The main purpose of the SAA approach and our modification is to 
reduce the large number of discretized PDE solves required for the inversion. 
In Table~\ref{Table: All the numbers for Exp1_newOpt}, we give a comparison of the total number of PDE solves for Example 1. Our approach drastically reduces the number of large-scale linear systems that needs to be solved. 
Additionally, it substantially improves the reconstruction results of the SAA approach.

\begin{table}
\begin{tabular}{c | c  c  c  c  c}
     									& \begin{tabular}{@{}c@{}} Iteration \\Number \end{tabular} & \begin{tabular}{@{}c@{}} Function \\Evaluations \end{tabular} &  \begin{tabular}{@{}c@{}} Jacobian \\Evaluations \end{tabular} & \begin{tabular}{@{}c@{}} Total PDE \\ Solves \end{tabular}  &  Tol\\ \hline
SAA$^\star$  (intermediate tol)		   & 10   & 11     & 6      & 170      & 	$\delta$		\\
1 Opt simult src/det       		& 18   & 19     & 10     & 524    &		 $\delta^2$	\\
2 Opt simult srcs/dets    	  & 18   & 19     & 10      &524  &		 $\delta^2$		\\
3 Opt simult srcs/dets        & 16   & 17     & 8      &484   &		 $\delta^2$	\\
All srcs/All dets 	&71    & 72      & 47         & 3808   &     $\delta^2$	\\
SAA$^{\star\star}$			&(32)    & (33)      & (19)         &  (520)  &     $\delta^2$	\\
SAA	$^{\star\star\star}$						&(92)    & (93)      & (67)         &  (1700)  &     $ \delta^2$	\\
\end{tabular}
\caption{ Example 1 Results. The total number of iterations, function evaluations, Jacobian evaluations and PDE solves required on average for 50 trials to reach the stopping criterion, 
$\| \Vr(\Vp) \|_2^2 = \delta^2$. \\
$^\star$The first row gives the cost to reach the intermediate tolerance for the SAA approach, 
$\| \Vr(\Vp) \|_2^2 = \delta$.
$^{\star\star}$	Since the SAA estimate becomes biased and underestimates the objective function, the algorithm stops prematurely.  $^{\star\star\star}$The SAA approach measuring the convergence with the true objective function. Parentheses indicate that the SAA approach does not reach the tolerance. \\}
\label{Table: All the numbers for Exp1_newOpt}
\end{table}

\begin{table}
\centering
\begin{tabular}{|l|l|l|l|l|l|}
\hline
\multicolumn{3}{|c|}{SAA Approach}                               & \multicolumn{3}{c|}{Rand $\&$ Optimized Simult Src/Det}                    \\ \cline{1-6}
                            Iter & True $\|\res \|_2^2$  ($\delta^2$) & Estimated  $\|\res \|_2^2$  ($\delta^2$) & Iter & True $\|\res \|_2^2$  ($\delta^2$)& Estimated $\|\res\|_2^2$  ($\delta^2$) \\ \hline     
 1    &118940         	 & 38820          &  1-5       		   &   (SAA)$^\star$       &  (SAA)$^\star$   \\
                             6    & 1192.5            & 391.15           &  6   				   & 1194.9  					  & 1197.3			\\                                                                                
                            11   &56.550            & 22.575           & 13   			    	& 118.73 					  & 118.56            \\
                            15   & 3.748             &0.8650         & 16   						& 19.894 						& 19.929         \\
                           (99)   & $ -$              & $ -$ 		& 22  						&0.8403 						&  0.8389
                                                                  \\ \hline
\end{tabular}
\caption{Subset of results for Example 1. The comparison of the true objective function $\|\res \|_2^2$ and its SAA estimate relative to the stopping criterion ($\delta^2$) for selected iterations. For the SAA approach, the estimated residual is obtained with 10 random simultaneous sources and detectors. Parentheses indicate that the SAA approach does not reach the tolerance. The estimated residual norms using  
random and optimized simultaneous sources and detectors are obtained using 3 optimized simultaneous sources and detectors. (SAA)$^\star$ indicates that we initially use the SAA approach up to the intermediate tolerance.}
\label{Table: Showing biased/unbiased_newOpt}
\end{table}

\el
\textbf{3D Experiment.} We use a $32\times 32\times 32$ grid, which gives $32,768$ unknowns 
in the discretized PDE (\ref{eq: PDE}). The model has 225 sources at the top and 225 detectors on the bottom, and we use only the zero frequency. In the PaLS approach, we use 27 CSRBFs, which leads to 135 parameters (five per 3D basis function) for the nonlinear optimization. The absorption image using the initial set of parameters is given in Figure~\ref{Fig: Initial Conf}b where 13 basis functions have a positive expansion coefficient (visible as high absorption regions) and 14 basis functions have a negative expansion coefficient (invisible). In our approach, we use only 12 random simultaneous  sources and detectors.

\begin{figure}
    \centering
       \begin{subfigure}{0.32\textwidth}
        \centering
        \includegraphics[ width=0.95\textwidth]{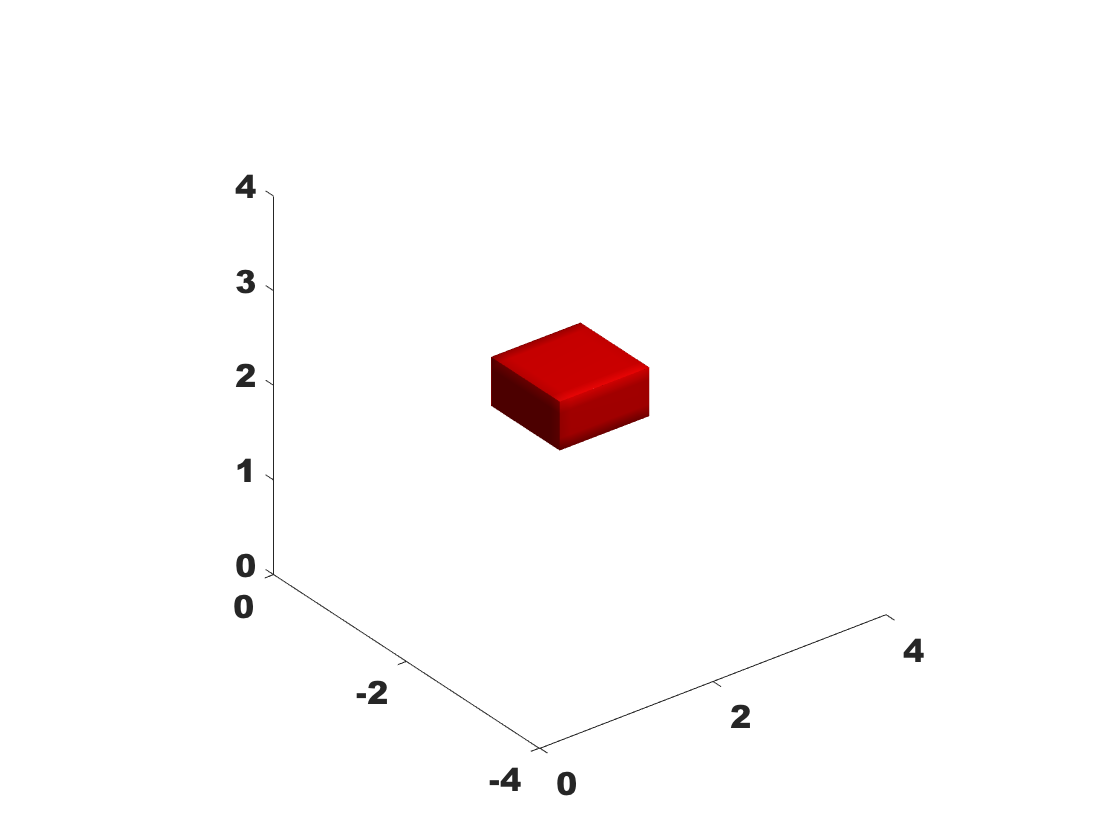}
        \subcaption{ }
    \end{subfigure}
        \centering
       \begin{subfigure}{0.32\textwidth}
        \centering
        \includegraphics[width=0.95\textwidth]{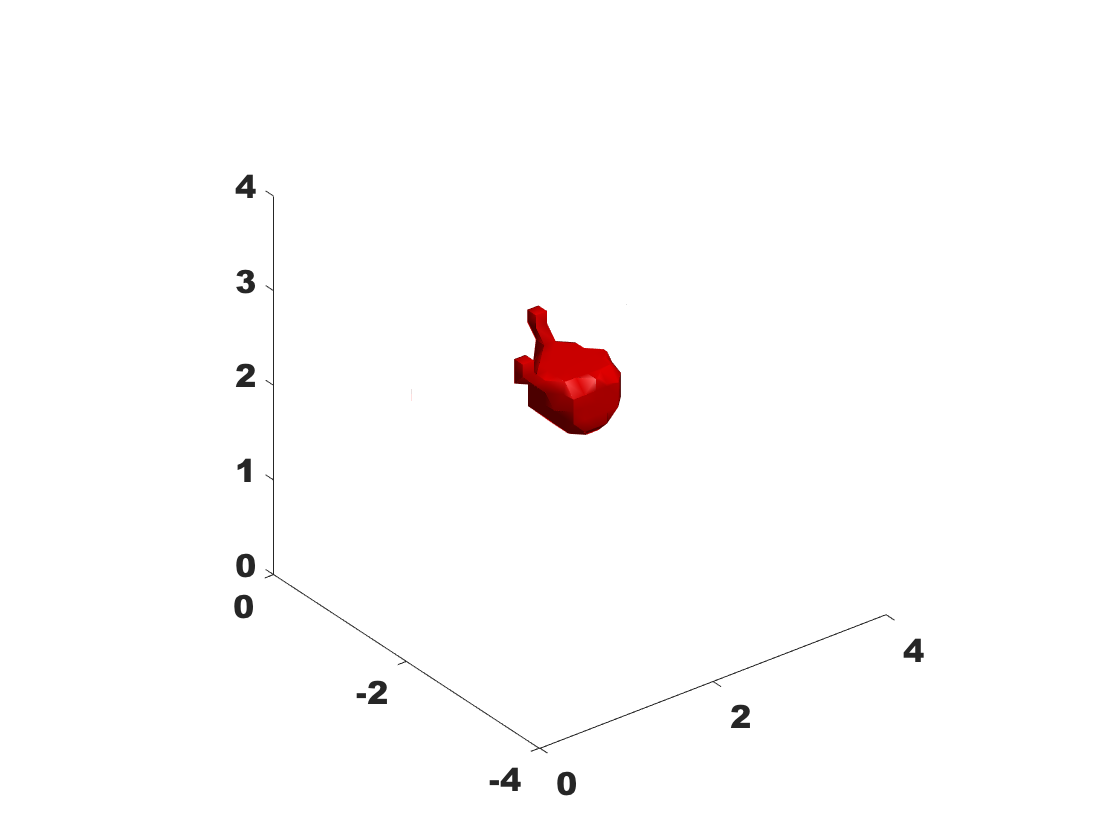}
        \subcaption{ }
    \end{subfigure}
    \centering
          \begin{subfigure}{0.32\textwidth}
        \centering
        \includegraphics[width=0.95\textwidth]{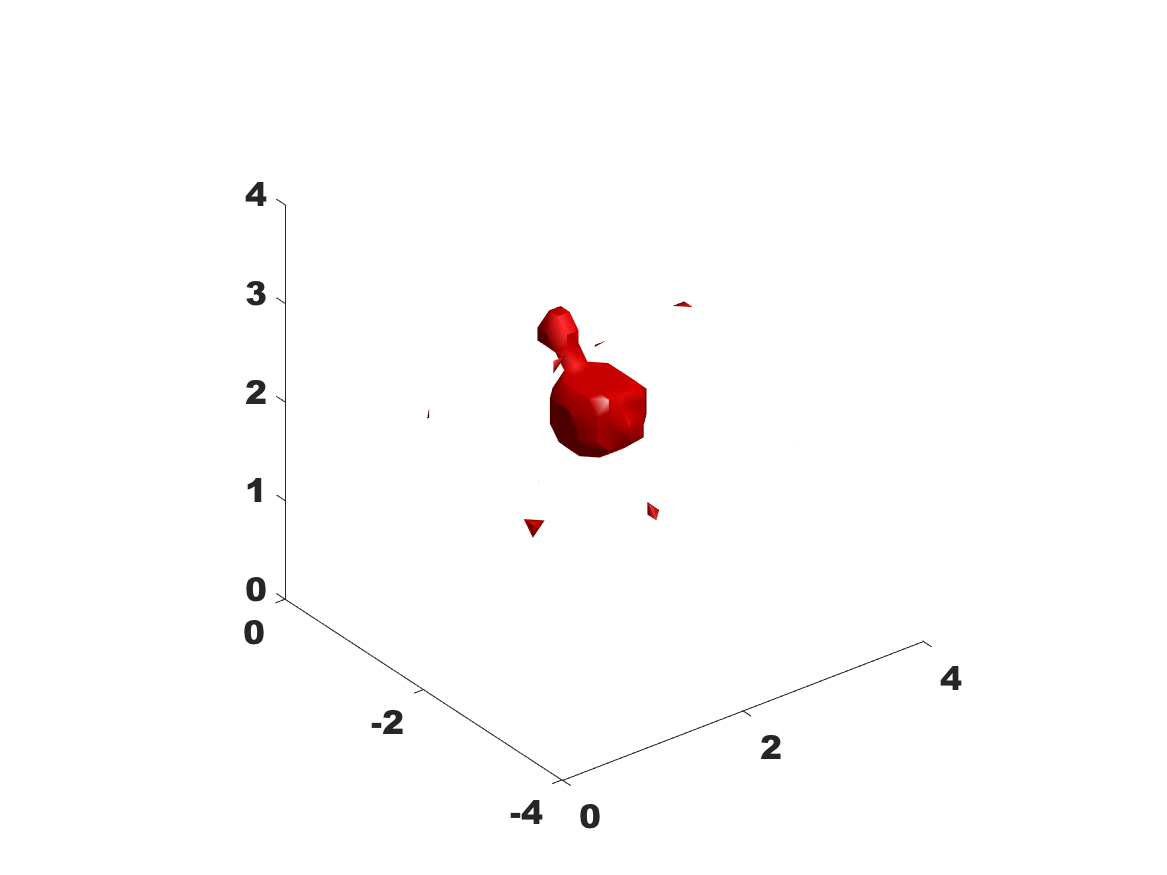}
        \subcaption{ }
    \end{subfigure}
    \centering
       \begin{subfigure}{0.32\textwidth}
        \centering
        \includegraphics[width=0.95\textwidth]{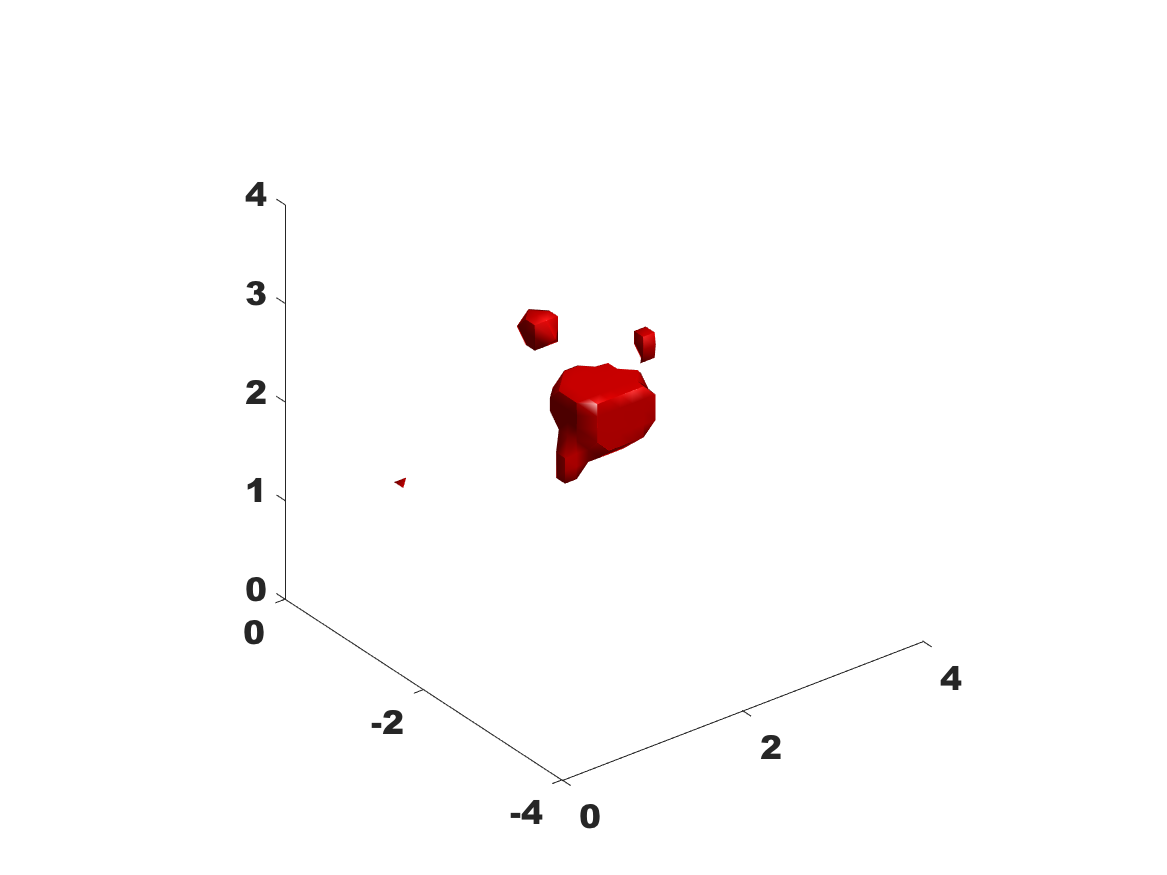}
        \subcaption{ }
    \end{subfigure}
       \centering
          \begin{subfigure}{0.32\textwidth}
        \centering
        \includegraphics[width=0.95\textwidth]{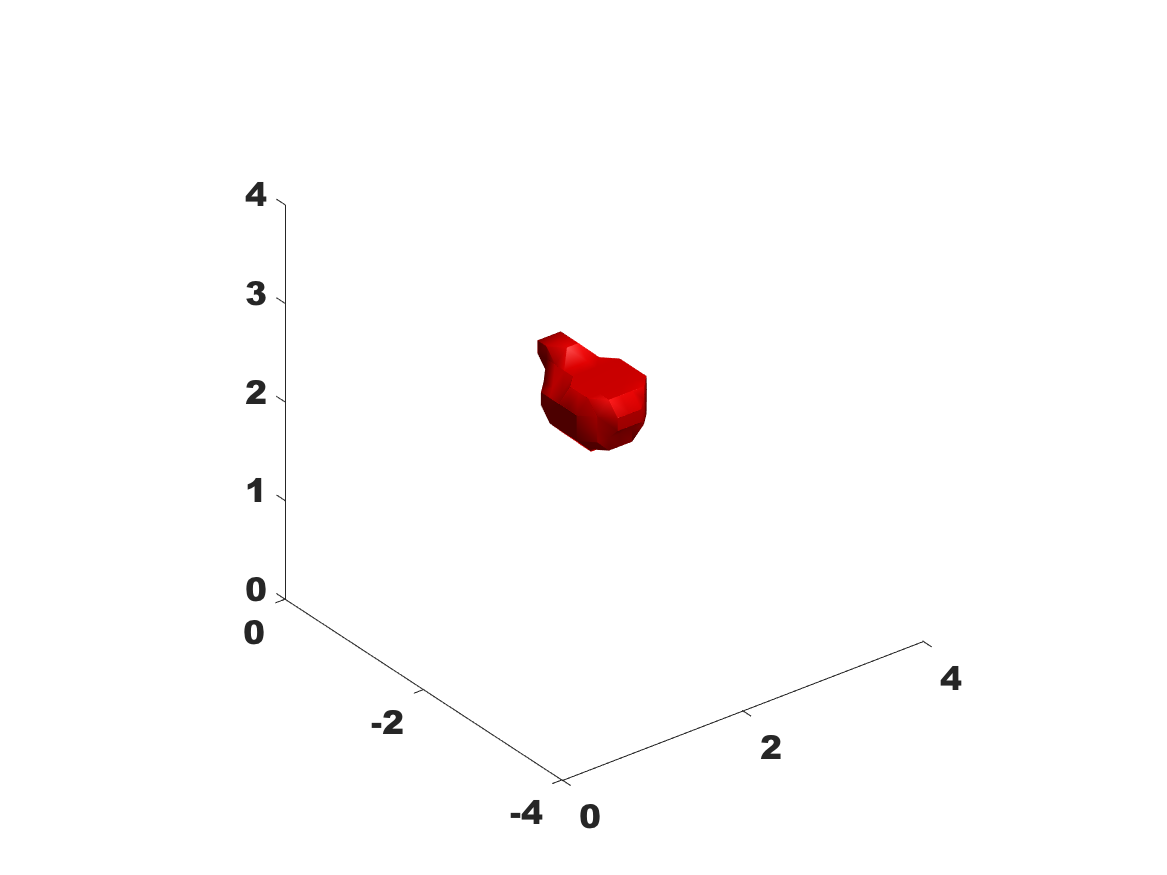}
        \subcaption{ }
    \end{subfigure}
    \centering
       \begin{subfigure}{0.32\textwidth}
        \centering
        \includegraphics[width=0.95\textwidth]{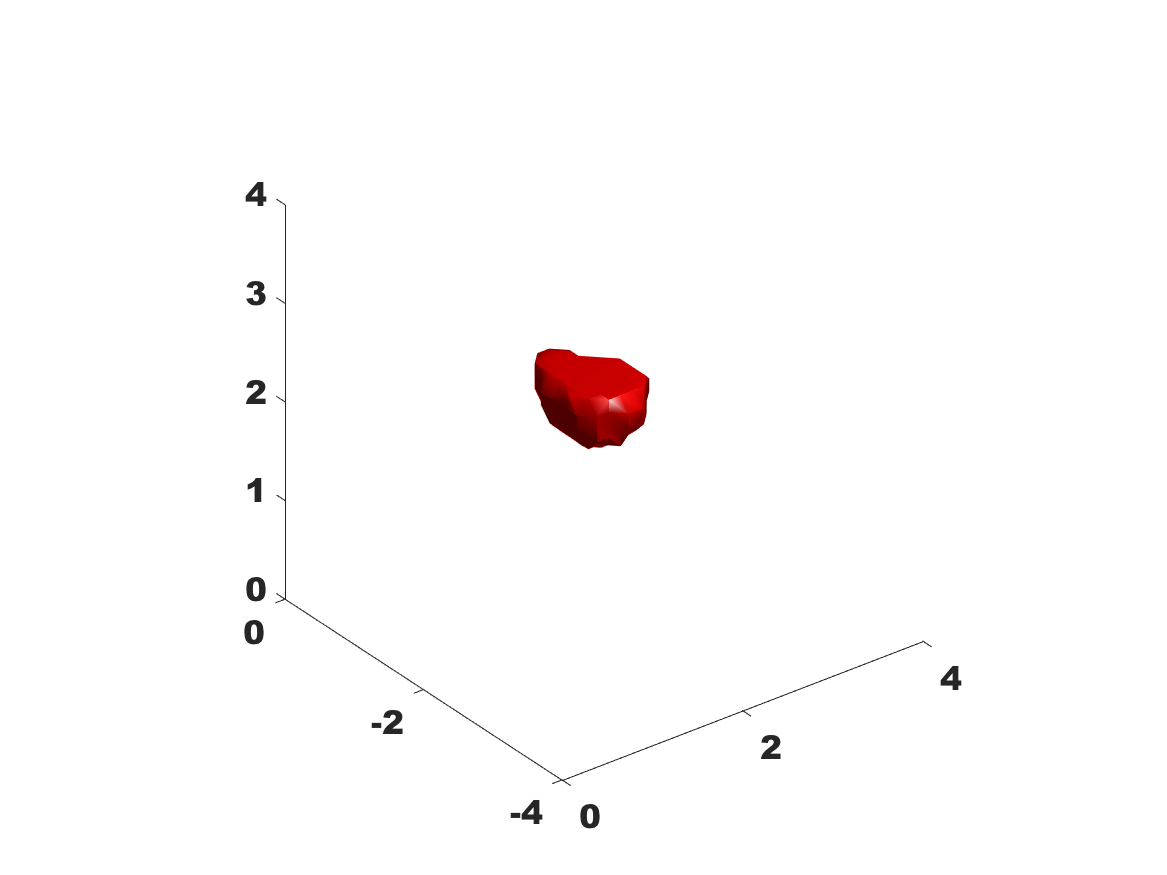}
        \subcaption{ }
    \end{subfigure}
\caption{Results for Example 2. Reconstruction of a test anomaly on a 
$32\times 32 \times 32$ mesh with 225 sources and detectors, 27 basis functions, and using only the zero frequency. The SAA approach uses 15 random simultaneous sources and detectors.\\
(a) True shape of the anomaly. (b) Reconstruction using all sources and all detectors. (c) Reconstruction 
using the SAA approach at the intermediate tolerance.
(d) Reconstruction using the SAA approach after the maximum iterations. (e) Reconstruction with SAA and 2 optimized simultaneous sources and detectors. (f) Reconstruction with SAA and 4 optimized simultaneous sources and detectors.}
\label{Fig:Recons All src det vs SAA vs Repl Exp3_newOpt}
\end{figure}

\textbf{Example 2.} The true absorption image for Example 2 is given in Figure~\ref{Fig:Recons All src det vs SAA vs Repl Exp3_newOpt}a. The reconstruction using all sources and detectors is given in Figure~\ref{Fig:Recons All src det vs SAA vs Repl Exp3_newOpt}b. Figure~\ref{Fig:Recons All src det vs SAA vs Repl Exp3_newOpt}c shows that SAA approach gives a good localization of the anomaly  at the intermediate tolerance. However, little further improvement appears using the SAA approach, even after the maximum iterations; see Figure~\ref{Fig:Recons All src det vs SAA vs Repl Exp3_newOpt}d. In Figure~\ref{Fig:Recons All src det vs SAA vs Repl Exp3_newOpt}e-g, we show the reconstruction results when combining random and optimized simultaneous sources and detectors.

The straightforward inversion using all sources and detectors requires $9,225$ large linear solves. Table \ref{Table: All the numbers for Exp3_newOpt} shows that our approach reduces the number of large linear solves by \textit{about a factor 12} compared with using all sources and detectors, while approximating the original shape well. Clearly, there is a large
improvement to be gained by using a small number of optimized simultaneous sources and detectors. For larger problems with many sources and detectors and using multiple frequencies, we expect much larger gains.

Overall, our approach improves the rate of convergence of the optimization and reduces the number of large-scale linear systems solves. Moreover, combining random and optimized simultaneous sources and detectors improves the quality of the inverse solution.

\begin{table}
\begin{tabular}{c | c  c  c  c  c}
     									& \begin{tabular}{@{}c@{}} Iteration \\Number \end{tabular} & \begin{tabular}{@{}c@{}} Function \\Evaluations \end{tabular} &  \begin{tabular}{@{}c@{}} Jacobian \\Evaluations \end{tabular} & \begin{tabular}{@{}c@{}} Total PDE \\ Solves \end{tabular}  &  Tol\\ \hline
SAA$^\star$  (intermediate tol)		   & 4   & 5    & 6      & 132      & 	$\delta$		\\
2 Opt simult srcs/dets      					& 12   & 13     & 6     & 726    &		 $\delta^2$	\\
4 Opt simult srcs/dets    					 & 9   & 10     & 3      &762  &		 $\delta^2$		\\
All srcs/All dets 								&25    & 26      & 15      & 9225   &     $\delta^2$	\\
SAA	$^{\star\star}$				&(99)    & (100)      & (67)         &  (2505)  &     $ \delta^2$	\\
\end{tabular}
\caption{ Example 2 Results. The total number of iterations, function evaluations, Jacobian evaluations and PDE solves required to reach the stopping criterion, $\| \Vr(\Vp) \|^2_2 = \delta^2$.\\
$^\star$The first row gives the costs to reach the intermediate tolerance for the SAA approach, 
$\| \Vr(\Vp) \|_2^2 = \delta$. $^{\star\star}$The SAA approach measuring the convergence with the true objective function. Parentheses indicate that the SAA approach does not reach the tolerance.}
\label{Table: All the numbers for Exp3_newOpt}
\end{table}
\section{Conclusions and Future Work}\label{Sec:SimultOpt_con}
We use the SAA approach to estimate the objective function, the Jacobian, and 
the gradient using only a few simultaneous random sources and detectors in 
DOT problems. While this approach is reasonably effective for the application 
in \citep{HabeChunHerr2012}, it does not work quite that well for DOT. 
Since convergence slows down in later iterations before the noise level is 
reached, and the standard SAA approach regularly
does not converge to the noise level, we propose using optimized 
simultaneous sources and detectors. 
With the addition of optimized directions, we observe faster convergence, 
good quality reconstructions, and robustness. This technique could be quite useful 
in other applications as well.

Several further improvements should be considered in the future.
In particular, approximating the, typically low rank, Jacobian at low cost but sufficiently
accurately to compute effective optimized simultaneous sources and detectors would
lead to a further substantial reduction in the number of large linear
solves. Potentially, computing such a low rank approximation can be combined
more efficiently with the tensor form of the Jacobian.

Although our approach has proved successful experimentally, we aim to understand the underlying theory better. In the future, we plan to analyze, more fundamentally, what are the most effective simultaneous sources and detectors for fast convergence of the inverse problem: randomized, optimized (and in what sense), and their combination.
An alternative approach to improve convergence, studied in multiple papers 
\cite{FarbodDoelAscher2,ByrdChinNocedalWu, BollaByrdNocedal, FarbodMahoney1, FarbodMahoney2},
is to slowly increase the sample size as the optimization progresses or dynamically choose
the sample size. In future work, we plan to compare these approaches with the approach proposed in this paper.
We also plan to test and evaluate SA approaches for inversion in DOT.

We intend to update the TREGS algorithm and study how small we can make the number of simultaneous sources and detectors (random and optimized) and still obtain good solutions and fast convergence. Moreover, finding more appropriate stopping criteria for the randomized approach may also improve our results.

As shown in \cite{StuGuKilChatBeatOCon}, using parametrized interpolatory model 
reduction can also reduce the cost of inversion for DOT. 
We plan to combine model reduction with the randomized approach.

\section{Acknowledgements}
We sincerely thank the reviewers for many helpful suggestions that 
substantially improved this paper.


\bibliographystyle{abbrv}
\bibliography{ms}

\end{document}